\documentclass[12pt,titlepage,english]{article}
\usepackage{amsmath}
\usepackage{amsfonts, amsthm}
\usepackage{amssymb}
\usepackage{amsbsy}
\usepackage{amstext}
\usepackage{booktabs}
\usepackage{caption}
\usepackage[usenames, dvipsnames]{color}
\usepackage[letterpaper]{geometry}
\usepackage{graphicx}
\usepackage{hyperref}
\usepackage[latin9]{inputenc}
\usepackage{lscape}
\usepackage{mathtools}
\usepackage{natbib}
\usepackage{setspace}
\usepackage{stmaryrd}
\usepackage{subcaption}
\usepackage{tabularx}
\usepackage[normalem]{ulem}
\usepackage{verbatim}

\counterwithin*{equation}{section}

\newtheorem{theorem}{Theorem}

\newtheorem{definition}{Definition}

\newtheorem{lemma}{Lemma}

\newtheorem{proposition}{Proposition}

\graphicspath{{figures/}{figures/homoskedastic_reduced/text/}{figures/hac/text/}}

\input{tcilatex}

\begin{document}

\title{Optimal Invariant Tests in an Instrumental Variables Regression With
Heteroskedastic and Autocorrelated Errors \thanks{%
Preliminary results were presented at seminars organized by BU, Brown,
Caltech, Harvard-MIT, PUC-Rio, University of California (Berkeley, Davis,
Irvine, Los Angeles, Santa Barbara, and Santa Cruz campuses), UCL, USC, and
Yale, at the FGV Data Science workshop, and at conferences organized by
CIREq (in honor of Jean-Marie Dufour), Harvard University (in honor of Gary
Chamberlain), Oxford University (New Approaches to the Identification of
Macroeconomic Models), the Tinbergen Institute (Inference Issues in
Econometrics), and Vanderbilt (Identification in Econometrics. This study
was financed in part by the Coordenação de Aperfeiçoamento de Pessoal de N\'
ivel Superior - Brazil (CAPES) - Finance Code 001. This research was also
supported in part by the University of Pittsburgh Center for Research
Computing through the resources provided.}}
\author{Marcelo J. Moreira \\
\emph{FGV}\\
\and Mahrad Sharifvaghefi \\
\emph{University of Pittsburgh} \and Geert Ridder \\
\emph{USC}}
\date{\bigskip \today\\
}
\maketitle

\begin{abstract}
This paper uses model symmetries in the instrumental variable (IV)
regression to derive an invariant test for the causal structural parameter.
Contrary to popular belief, we show that there exist model symmetries when
equation errors are heteroskedastic and autocorrelated (HAC). Our theory is
consistent with existing results for the homoskedastic model (%
\citet{AndrewsMoreiraStock06} and \citet{Chamberlain07}). We use these
symmetries to propose the conditional integrated likelihood (CIL) test for
the causality parameter in the over-identified model. Theoretical and
numerical findings show that the CIL\ test performs well compared to other
tests in terms of power and implementation. We recommend that practitioners
use the Anderson-Rubin (AR) test in the just-identified model, and the CIL
test in the over-identified model.
\end{abstract}



\setcounter{equation}{0}

\section{Introduction \label{Sec Introduction}}

In a regression model, the explanatory variable can be correlated with the
error due to omitted variables. To solve this endogeneity problem,
practitioners often look for instrumental variables (IVs). The instruments
are valid if they are correlated with the endogenous variable but
uncorrelated with the error. The instruments are said to be weak when their
correlation with the endogenous explanatory variable is small. Under weak
identification, standard estimators may be far from the true causality
parameter, and commonly-used tests do not have correct size. Searching for
valid IVs can, unfortunately, narrow down the choices to only weak
instruments. Furthermore, techniques proposed to mitigate these problems can
themselves have limitations. \citet{CruzMoreira05} show that the
second-order bias for the two-stage least squares (2SLS) estimator is
unreliable under weak identification. \citet{LeeMcCraryMoreiraPorter20}
point out that the standard F>10 rule for the t-ratio leads to important
size distortions in practice, even in the just-identified model. They
propose a novel tF procedure if practitioners wish to use the F statistic
combined with the t-ratio.

With cross-sectional data, the errors in the IV model can be
heteroskedastic. With time-series or panel data, errors can also be
autocorrelated. For these more complex data generating processes (DGPs) for
the errors, the asymptotic variance matrix of sample IV moments can be quite
different from the one obtained from serially uncorrelated and homoskedastic
errors. Consistent estimators for this variance are readily available: see %
\citet{NeweyWest87} and \citet{Andrews91}.

\citet{AndrewsMoreiraStock06} (abbreviated as AMS06 hereinafter) and %
\citet{Chamberlain07} show that symmetries exist in the IV model with
homoskedastic errors when the variance is fixed. Up until the emergence of
this project (\citet{MoreiraRidder17}), it was widely accepted that there
are no symmetries in the IV model with heteroskedastic and/or autocorrelated
(HAC) errors. Indeed, at first sight, invariance does not seem applicable to
the HAC-IV model. We argue that this view is incorrect. The HAC-IV has model
symmetries if the HAC variance matrix is assumed to be known but not fixed,
an important distinction from the method used by AMS06, \citet{Chamberlain07}%
, and related papers, for homoskedastic and uncorrelated errors. We find the
largest affine group which preserves the null hypothesis for the causality
parameter. This allows us to find weights for the novel conditional
integrated likelihood (CIL)\ test. This test is invariant and can be
interpreted as the limit of a sequence of conditional weighted-average power
(WAP)\ tests. \citet{AndrewsMikusheva20} provide a general framework for
decision rules in GMM. They do not propose a two-sided test, which is our
main goal in this paper. Unlike the CIL\ test, other (limits of) conditional
WAP\ tests can be severely biased, as the critique by %
\citet{MoreiraMoreira19} asserts.

In the just-identified model, the AR test is optimal within the classes of
either unbiased or invariant tests, assuming the reduced-form variance
mentioned above is known; see \citet{Moreira02, Moreira09a}, AMS06, and %
\citet{MoreiraMoreira19}. In the over-identified model, the AR test is not
efficient under the usual asymptotic theory. Several proposed tests are
asymptotically optimal under standard asymptotics. \citet{MoreiraRidder20}
show that the Lagrange Multiplier/score (LM) and the conditional
quasi-likelihood ratio (CQLR) tests can suffer severe power deficiencies
when distinguishing the null from the alternative hypothesis should be easy.
The weighted-average strongly unbiased (SU) tests are not invariant for
arbitrary weight choices. Furthermore, implementing the SU tests requires
linear programming. Although algorithms are readily available, it requires
the calculation of a density ratio. If not dealt with properly, the
computation can exceed the numerical accuracy of computer packages. This
leaves, as the main contender, the true conditional likelihood ratio (CLR)
test; see \citet{AndrewsMikusheva16} and \citet{MoreiraMoreira19}. The CLR
test does not have a closed-form expression with HAC errors, and requires
numerical optimization. We show some important limitations to the
implementation of the CLR test. We prove here that we can compactify the
parameter space for the optimization. This is important, and avoids some
natural pitfalls if the parameter space is unbounded. However, we find that
the number of initial points needed for the optimization algorithm depends
on the errors' DGPs and on the instrument strength. In practice, we document
the need to include several initial points when errors are HAC. Worse yet,
the optimization can be even more troublesome when computing the conditional
quantiles. This happens because the model is misspecified for DGPs under the
alternative hypothesis when we simulate these null conditional quantiles.

We then compare power between the AR, CLR, and CIL\ tests. For homoskedastic
errors, the CLR test simplifies to the CQLR test, which has a closed-form
expression. \citet{AndrewsMoreiraStock04} use these symmetries to choose
weights for a conditional WAP test. Power gains beyond the CQLR are small,
as the latter performs near a two-sided power envelope for invariant tests
with homoskedastic errors. It is, however, reassuring that the CIL test\
performs on an equal footing with the CQLR test. For HAC errors, the IV
model is much more complex than the simple homoskedastic model. Even
reducing the data using invariance, several parameters can affect the
performance of AR, LM, CQLR, CLR, CIL, and any other invariant tests. We
choose the same designs as \citet{MoreiraMoreira19} and %
\citet{MoreiraRidder20}, to forestall any criticism that we may be selecting
parameter combinations which favor the CIL test. To bypass the
aforementioned numerical problems for the CLR test, we choose to implement
an infeasible version of the CLR test, in case better optimization methods
are found in the future. This implementation selects the unknown value of
the structural parameter as one of the initial points in the likelihood
optimization. Overall, the CIL test outperforms the AR and CLR tests, with
significant power gains for several of these designs.

This paper is organized as follows. Section \ref{Sec IV Model and Statistics}
introduces the IV model, and describes the family of similar tests robust to
heteroskedastic-autocorrelated errors. Section \ref{Sec Invariance and CIL}
shows model symmetries when the asymptotic variance can change with data
transformations. We present different representations of the CIL\ test. One
of them is important to show that this test is invariant, as discussed
later. The other expression is useful to derive confidence sets based on the
CIL\ test. Section \ref{Sec Numerical} shows that the CIL\ test can have
very good power (the supplement provides further evidence in favor of the
CIL\ test). The more technical details behind model invariance are left to
the end of the paper. Section \ref{Sec Kronecker} shows that the theory of
AMS06 is a special case of ours when the variance has a Kronecker product
form. Section \ref{Sec Invariant Tests} derives the theory of conditional
invariant tests. It shows that the AR, CLR, and Lagrange Multiplier/score
(LM) tests are also invariant. Section \ref{Sec Conclusion} discusses the
next steps in this research agenda and highlights the methodological
importance of distinguishing between parameters being known and being fixed.
The online appendix provides the proofs for our theory.

\section{The IV Model and Statistics \label{Sec IV Model and Statistics}}

\subsection{The HAC-IV model \label{Sec HAC-IV}}

Consider the following structural equation for the $i$-th observation of the
variable $y_{1}$: 
\begin{equation}
y_{1i}=y_{2i}\beta +x_{i}^{\prime }\gamma _{1}+u_{i},\text{ for }%
i=1,2,\cdots ,n,  \label{dgp_y_1i}
\end{equation}%
where $y_{2i}$ is an endogenous random variable with corresponding
coefficient $\beta \in \mathbb{R}$, $x_{i}=(x_{1i},x_{2i},\cdots
,x_{pi})^{\prime }\in \mathbb{R}^{P}$ is a fixed vector of exogenous control
variables with corresponding vector of coefficients $\gamma _{1}=(\gamma
_{1}^{\ast },\gamma _{2}^{\ast },\cdots ,\gamma _{p}^{\ast })^{\prime }\in 
\mathbb{R}^{p}$, and $u_{i}$ is an error term. We also consider the
following reduced-form equation for the endogenous explanatory random
variable: 
\begin{equation}
y_{2i}=\tilde{z}_{i}^{\prime }\pi +x_{i}^{\prime }\xi _{1}+v_{2i},\text{ for 
}i=1,2,\cdots ,n,  \label{dgp_y_2i}
\end{equation}%
where $\tilde{z}_{i}=(\tilde{z}_{1i},\tilde{z}_{2i},\cdots ,\tilde{z}%
_{ki})^{\prime }\in \mathbb{R}^{k}$ is a fixed vector of instrumental
variables (IVs) with corresponding coefficients $\pi =(\pi _{1},\pi
_{2},\cdots ,\pi _{k})^{\prime }\in \mathbb{R}^{k}$, $\xi _{1}=(\xi _{1},\xi
_{2},\cdots ,\xi _{p})\in \mathbb{R}^{p}$, and an error term $v_{2i}$. It
may be possible that $\mathbb{E}(\nu _{2i}u_{i})\neq 0$, so that $y_{2}$ is
an endogenous random variable. Equations (\ref{dgp_y_1i}) and (\ref{dgp_y_2i}%
) can be presented in the following matrix format:%
\begin{eqnarray}
y_{1} &=&y_{2}\beta +X\gamma _{1}+u  \label{(Struc Eqns 1)} \\
y_{2} &=&\widetilde{Z}\pi +X\xi _{1}+v_{2},  \notag
\end{eqnarray}%
where $y_{1}=(y_{11},y_{12},\cdots ,y_{1n})^{\prime }\in \mathbb{R}^{n}$, $%
y_{2}=(y_{21},y_{22},\cdots ,y_{2n})^{\prime }\in \mathbb{R}^{n},$ $%
X=(x_{1},x_{2},\cdots ,x_{n})^{\prime }\in \mathbb{R}^{n\times p}$, $%
\widetilde{Z}=(\tilde{z}_{1},\tilde{z}_{2},\cdots ,\tilde{z}_{n})^{\prime
}\in \mathbb{R}^{n\times k}$, $u=(u_{1},u_{2},\cdots ,u_{n})^{\prime }\in 
\mathbb{R}^{n}$, and $v_{2}=(v_{21},v_{22},\cdots ,v_{2n})^{\prime }\in 
\mathbb{R}^{n}$. We assume that the matrix $\overline{Z}=[\widetilde{Z}:X]$
has full column rank $k+p$.

Our focus is on testing the null hypothesis $H_{0}:\beta =\beta _{0}$
against the two-sided alternative hypothesis $H_{1}:\beta \neq \beta _{0}$.
It is convenient to transform the IV matrix $\widetilde{Z}$ to $Z$ which is
orthogonal to matrix $X$, $Z^{\prime }X=0$. For a conformable matrix $A$, we
define $N_{A}=A(A^{\prime }A)^{-1}A^{\prime }$ and $M_{A}=I-N_{A}$. We write 
\begin{equation}
y_{2}=Z\pi +X\xi +v_{2},  \label{(Orth Red Form Eqn)}
\end{equation}%
where $Z=M_{X}\widetilde{Z}\text{ and }\xi =\xi _{1}+(X^{\prime
}X)^{-1}X^{\prime }\widetilde{Z}\pi $. By substituting $y_{2}$ from the
reduced-form equation (\ref{(Orth Red Form Eqn)}) to the structural equation
(\ref{(Struc Eqns 1)}), we have 
\begin{equation}
y_{1}=Z\pi \beta +X\gamma +v_{1},  \label{reduced from y1}
\end{equation}%
where $\gamma =\gamma _{1}+\xi \beta $ and $v_{1}=u+v_{2}\beta $. The
reduced-form equations (\ref{(Orth Red Form Eqn)}) and (\ref{reduced from y1}%
) can be written in the following matrix notation: 
\begin{equation}
Y=Z\pi a^{\prime }+X\eta +V,  \label{(Matrix Red Form)}
\end{equation}%
where $Y=[y_{1}:y_{2}]\in \mathbb{R}^{n\times 2}$, $V=[v_{1}:v_{2}]\in 
\mathbb{R}^{n\times 2}$, $a=(\beta ,1)^{\prime },\text{ and }\eta =[\gamma
:\xi ]\in \mathbb{R}^{p\times 2}$.

\Citet{MoreiraMoreira19} consider $R\equiv \left( Z^{\prime }Z\right)
^{-1/2}Z^{\prime }Y\in \mathbb{R}^{k\times 2}$. Because the transformed IV
matrix $Z$ is orthogonal to $X$, we have%
\begin{equation}
R=\mu a^{\prime }+\widetilde{V},
\end{equation}%
where $\widetilde{V}=\left( Z^{\prime }Z\right) ^{-1/2}Z^{\prime }V$ and $%
\mu =\left( Z^{\prime }Z\right) ^{1/2}\pi $. Commonly-used estimators and
tests depend on the data through $R$ and estimators $\widehat{\Sigma }_{n}$
for the variance $\Sigma _{n}$ of $vec(\widetilde{V})$. For example,
consider the t-statistic based on the 2SLS estimator%
\begin{equation}
\widehat{\beta }=\left( y_{2}^{\prime }N_{Z}y_{2}\right) ^{-1}y_{2}^{\prime
}N_{Z}y_{1}.
\end{equation}%
The estimator is clearly a ratio of quadratic forms of $R$. In our notation,
the t-statistic (also known as the Wald statistic) is%
\begin{equation}
\widehat{W}_{n}=\frac{\widehat{\beta }-\beta _{0}}{\widehat{\sigma }_{\beta
,n}}\text{ for }\widehat{\sigma }_{\beta ,n}^{2}=\frac{R_{2}^{\prime }(%
\widehat{b}^{\prime }\otimes I_{k})\widehat{\Sigma }_{n}(\widehat{b}\otimes
I_{k})R_{2}}{\left( R_{2}^{\prime }R_{2}\right) ^{2}}\text{,}
\end{equation}%
where $\widehat{b}=\left( 1,-\widehat{\beta }\right) ^{\prime }$ and $R_{2}$
is the second column of $R$.

The two-sided t-test rejects the null when $\left\vert W\right\vert $ is
larger than the $1-\alpha $ quantile of a standard normal distribution. For
this critical value to be reliable, the t-statistic needs to be
approximately normally distributed. This happens when the number of
observations $n$ increases and the IVs are strong. In applied work, however,
it can be difficult to find variables that are also uncorrelated with the
error terms of the structural equation (\ref{dgp_y_1i}). In practice, the
search for valid IVs may lead to choices which are weakly correlated with
the endogenous explanatory variable $y_{2}$. As a result, the null rejection
probability for the t-test can be sensitive to the quality of the
instruments; see \Citet{NelsonStartz90b}, \Citet{Dufour97}, and %
\Citet{StaigerStock97}. In particular, the null rejection probability can be
much larger than the usual nominal level. This problem spurs us to develop
similar tests which, by construction, have null rejection probability equal
to nominal level $\alpha $, no matter how weak the IVs are.

\subsection{Similar Tests \label{Sec Statistics and Tests}}

For simplicity, we start by assuming that $vec(\widetilde{V})$ is normally
distributed with zero mean and $\Sigma $ is known. The online appendix
relaxes this assumption, at the cost of asymptotic approximations. For
example, the t-statistic for known $\Sigma $ (to streamline notation, we
omit the subscript $n$ from $\Sigma _{n}$) would be%
\begin{equation}
W=\frac{\widehat{\beta }-\beta _{0}}{\widehat{\sigma }_{\beta }}\text{,
where }\widehat{\sigma }_{\beta }^{2}=\frac{R_{2}^{\prime }(\widehat{b}%
^{\prime }\otimes I_{k})\Sigma (\widehat{b}\otimes I_{k})R_{2}}{\left(
R_{2}^{\prime }R_{2}\right) ^{2}}.
\end{equation}%
For other test statistics, it is convenient to transform $R$ into the pair
of $k\times 1$ statistics, $S$ being pivotal and independent of the
statistic $T$. \Citet{MoreiraMoreira19} and \Citet{MoreiraRidder20} define%
\begin{align}
S& =\left[ \left( b_{0}^{\prime }\otimes I_{k}\right) \Sigma \left(
b_{0}\otimes I_{k}\right) \right] ^{-1/2}\left( b_{0}^{\prime }\otimes
I_{k}\right) vec\left( R\right) \text{ and}  \label{(Defns of S and T)} \\
T& =\left[ \left( a_{0}^{\prime }\otimes I_{k}\right) \Sigma {}^{-1}\left(
a_{0}\otimes I_{k}\right) \right] ^{-1/2}\left( a_{0}^{\prime }\otimes
I_{k}\right) \Sigma ^{-1}vec\left( R\right) \text{,}  \notag
\end{align}%
for $a_{0}=\left( \beta _{0},1\right) ^{\prime }$ and $b_{0}=\left( 1,-\beta
_{0}\right) ^{\prime }$. Their marginal distributions are given by%
\begin{eqnarray}
S &\sim &N\left( \left( \beta -\beta _{0}\right) C_{\beta _{0}}\mu
,I_{k}\right) \text{ and }T\sim N\left( D_{\beta }\mu ,I_{k}\right) \text{,
where}  \label{(Dist S and T)} \\
C_{\beta _{0}} &=&\left[ \left( b_{0}^{\prime }\otimes I_{k}\right) \Sigma
\left( b_{0}\otimes I_{k}\right) \right] ^{-1/2}\text{ and}  \notag \\
D_{\beta } &=&\left[ \left( a_{0}^{\prime }\otimes I_{k}\right) \Sigma
^{-1}\left( a_{0}\otimes I_{k}\right) \right] ^{-1/2}\left( a_{0}^{\prime
}\otimes I_{k}\right) \Sigma ^{-1}\left( a\otimes I_{k}\right) .  \notag
\end{eqnarray}

Examples of test statistics based on $S$ and $T$ are the Anderson-Rubin
(AR), the score or Lagrange multiplier (LM), and the quasi likelihood ratio
(QLR) statistics. \citet{AndersonRubin49} propose a pivotal test statistic.
In our model, the Anderson-Rubin statistic is given by 
\begin{equation}
AR=S^{\prime }S.  \label{(AR stat)}
\end{equation}%
\citet{MoreiraMoreira19} derive the $LM$ statistic under the same
distributional assumption that we make here. The two-sided $LM$ statistic is%
\begin{equation}
LM=S^{\prime }N_{C_{\beta _{0}}D_{\beta _{0}}^{-1}T}S.  \label{(LM stat)}
\end{equation}%
The $AR$ and $LM$ statistics have chi-square distributions with $k$ and one
degrees of freedom, respectively. The AR and LM\ tests reject the null when
their respective statistics are larger than their $1-\alpha $ chi-square
quantiles. By construction, both tests have correct size at level $\alpha $.

\citet{Kleibergen05}, among others, adapts the likelihood ratio statistic
for homoskedastic errors to HAC errors. The quasi-likelihood ratio statistic
is%
\begin{equation}
QLR=\frac{AR-r\left( T\right) +\sqrt{\left( AR-r\left( T\right) \right)
^{2}+4LM\cdot r\left( T\right) }}{2},  \label{(QLR stat)}
\end{equation}%
where $r\left( T\right) =T^{\prime }T$ . \citet{Andrews16} proposes tests
based on the following combination:%
\begin{equation}
LC=m\left( T\right) \cdot AR+\left( 1-m\left( T\right) \right) \cdot LM,
\end{equation}%
where $0\leq m\left( T\right) \leq 1$. Unlike the $AR$ and $LM$ statistics,
neither the $QLR$ nor the $LC$ statistics are pivotal. We follow %
\citet{MoreiraMoreira19} and reject the null hypothesis when the test
statistic $\psi $ is larger than $\kappa \left( t,\Sigma \right) $, which is
the null $1-\alpha $ quantile conditional on $T=t$. Writing a test statistic
as $\psi \left( S,T,\Sigma \right) ,$ we can compute the conditional
rejection probability under the null:%
\begin{equation}
P_{\beta _{0},\mu ,\Sigma }\left( \psi \left( S,T,\Sigma \right) \geq
x|T=t\right) =P_{\beta _{0},\mu ,\Sigma }\left( \psi \left( S,t,\Sigma
\right) \geq x\right) .
\end{equation}%
This probability does not depend on $\mu $ because the distribution of $S$
under the null is pivotal. By construction, the conditional null quantile
satisfies%
\begin{equation}
P_{\beta _{0},\mu ,\Sigma }\left( \psi \left( S,t,\Sigma \right) \geq \kappa
\left( t,\Sigma \right) \right) \equiv \alpha .
\end{equation}%
Consequently, the unconditional null rejection probability is $\alpha $,%
\begin{equation}
P_{\beta _{0},\mu ,\Sigma }\left( \psi \left( S,T,\Sigma \right) \geq \kappa
\left( T,\Sigma \right) \right) \equiv \alpha .
\end{equation}%
For example, the conditional test based on the $QLR$ statistic rejects the
null when this statistic is larger than its null conditional quantile. If
the statistic is pivotal, like the $AR$ and $LM$ statistics, the conditional
quantile $\kappa \left( t,\Sigma \right) $ collapses to the null
unconditional quantile.

The $QLR$ and $LC$ statistics depend on $S$ only through the $AR$ and $LM$
statistics. \citet{MoreiraRidder20} show that the statistic $S$ has useful
information beyond the Anderson-Rubin and score statistics when the
covariance matrix does not have a Kronecker product structure. For that
reason, we recommend the use of conditional tests based on either a
likelihood ratio statistic or a WAP statistic to be introduced here. These
tests take advantage of information beyond the Anderson-Rubin and score
statistics.

The likelihood ratio statistic based on $R$ is%
\begin{equation}
LR=\max_{a}vec(R)^{\prime }\Sigma ^{-1/2}N_{\Sigma ^{-1/2}(a\otimes
I_{k})}\Sigma ^{-1/2}vec(R)-T^{\prime }T,  \label{(LR stat)}
\end{equation}%
where $LR$ can be written in terms of the pivotal statistic $S$ and the
complete statistic $T$; see \citet{MoreiraMoreira19}. In the appendix, we
show that this statistic can be written as%
\begin{equation}
LR=b_{0}^{\prime }R^{\prime }\left[ \left( b_{0}^{\prime }\otimes
I_{k}\right) \Sigma \left( b_{0}\otimes I_{k}\right) \right]
^{-1}Rb_{0}-\min_{b}b^{\prime }R^{\prime }\left[ \left( b^{\prime }\otimes
I_{k}\right) \Sigma \left( b\otimes I_{k}\right) \right] ^{-1}Rb,
\label{(LR moment)}
\end{equation}%
where $b=\left( 1,-\beta \right) ^{\prime }$. Hence, $LR$ is associated to
the GMM objective function based on the moment $E\left( Z^{\prime }\left(
y_{1}-y_{2}\beta \right) \right) =0$ and the continuously-updating weighting
matrix; see \citet{AndrewsMikusheva16} for the general case. The $LR$
statistic does not have a closed-form solution and requires numerical
searching methods. We instead use invariance to find an integrated
likelihood test.

\section{Invariance and the CIL Test \label{Sec Invariance and CIL}}

Contrary to popular belief, the IV model with HAC errors presents symmetries.%
\footnote{%
An econometric model is a (parametric, semi-parametric, or non-parametric)
family $\mathcal{P}$ of probability measures $P$ for the data $Y$. Consider
the transformations on the data $g\circ Y$ given by a group $g\in \mathcal{G}
$. This action yields a transformation $g\circ P$ given by $g\circ P\left(
Y\in B\right) \equiv P\left( g\circ Y\in B\right) $ for any Borel set $B$.
The model is said to be symmetric when $g\circ P\in \mathcal{P}$ for every $%
g\in \mathcal{G}$ and $P\in \mathcal{P}$.} The theory developed for the IV
model thus far assumes the variance matrix is fixed. This assumption
prevents us from finding symmetries with more general error DGPs. In this
paper, we instead assume that the variance $\Sigma $ is known, but not fixed.

To explain the symmetries present in the IV model, first consider a simple
example, in which $Y_{i}\overset{iid}{\sim }N\left( \tau ,\sigma ^{2}\right) 
$, where $\sigma ^{2}$ is unknown. We want to test the null hypothesis $%
H_{0}:\tau =0$ against $H_{1}:\tau \neq 0$, treating $\sigma ^{2}$ as a
nuisance parameter. For any scalar $g\neq 0$, the transformed data $%
X_{i}=g\cdot Y_{i}$ has distribution $X_{i}\overset{iid}{\sim }N\left(
g\cdot \tau ,g^{2}\sigma ^{2}\right) $. This simple model is then symmetric
(or said to be preserved) for the multiplicative group $\mathcal{G}$. The
transformation preserves the null (and therefore, the alternative) because
the mean of $X_{i}$ is zero if and only if the mean of $Y_{i}$ is zero. The
sufficient statistic for $\left( \tau ,\sigma ^{2}\right) $ is the sample
mean $\overline{Y}_{n}$ and the variance estimator $S_{Y}^{2}=n^{-1}\sum
\left( Y_{i}-\overline{Y}_{n}\right) ^{2}$. The transformation above induces
a change in the space of sufficient statistics: the pair $\overline{Y}_{n}$
and $S_{Y}^{2}$ become $\overline{X}_{n}=g\cdot \overline{Y}_{n}$ and $%
S_{X}^{2}=g^{2}S_{Y}^{2}$, respectively. If these transformations preserve
the hypothesis-testing problem and the original data are supportive of a
hypothesis, the transformed data should be equally supportive of the same
hypothesis. This is called the \emph{invariance principle}. Therefore, the
test statistic should be the same whether computed from the original or from
the transformed data; in other words, the test has to be invariant. Any
invariant test can be written as a function of the largest invariant
statistic. In this example, the maximal invariant is then $\overline{X}%
_{n}^{2}/S_{X}^{2}=\overline{Y}_{n}^{2}/S_{Y}^{2}$. Its distribution depends
only on $\tau ^{2}/\sigma ^{2}$ and has a monotone likelihood ratio
property. As a result, the uniformly most powerful invariant (UMPI)\ test
rejects the null when $\overline{Y}_{n}^{2}/S_{Y}^{2}$ is sufficiently
large. We refer interested readers to \citet{Eaton89} and %
\citet{LehmannRomano05} for the theory of optimal tests.

Now, consider instead the case in which $\sigma ^{2}$ is known. The
multiplicative group does not preserve the model if we assume $\sigma ^{2}$
to be fixed. We would have to consider a much smaller group in which $g=\pm
1 $ only (this restriction is in perfect analogy to the sign group defined
by AMS06, as we shall see in Section \ref{Sec Kronecker}). However, this
transformation only reduces the sufficient statistic to the maximal
invariant $\overline{Y}_{n}^{2}$ and $S_{Y}^{2}$. How, then, can we use the
model symmetries to obtain a further reduction? One possibility is to
distinguish the assumption of a known variance from the assumption of a
fixed variance. The distinction hinges on whether we actually know $\sigma
^{2}$ and treat it as fixed, even after we transform the data. If an
outsider tells us the value of $\sigma ^{2}$, this person would give a
different answer if we asked what the variance is after multiplying the data
by a nonzero scalar. The person reports a known, but not fixed, variance. We
can still get an optimal test if we restrict ourselves to unbiased tests.
Because our simple model belongs to a one-parameter exponential family, we
automatically find that the uniformly most powerful unbiased (UMPU) test
rejects the null hypothesis for large values of $\overline{Y}%
_{n}{}^{2}/\sigma ^{2}$.

Instead, we can take the variance $\sigma ^{2}$ as both part of the data and
the parameter space. The sufficient statistic is now the pair $\overline{Y}%
_{n}$ and $\sigma ^{2}$, while the parameters are $\tau $ and also $\sigma
^{2}$. The same multiplicative group transforms the sufficient statistic to $%
\overline{X}_{n}=g\cdot \overline{Y}_{n}$ and $g^{2}\sigma ^{2}$, and
induces a change in the mean from $\tau $ to $g\cdot \tau $ and the variance
from $\sigma ^{2}$ to $g^{2}\sigma ^{2}$. The maximal invariant is then $%
\overline{X}_{n}^{2}/\sigma ^{2}=\overline{Y}_{n}^{2}/\sigma ^{2}$. This
statistic has a noncentral chi-square distribution, where the noncentrality
parameter $\tau ^{2}/\sigma ^{2}$ is zero if and only if the null hypothesis
is true. Because this distribution also has a monotonic likelihood ratio
property, we again obtain a UMPI test that rejects the null hypothesis if $%
\overline{Y}_{n}^{2}/\sigma ^{2}$ is large.

In this simple canonical model, the UMPU and UMPI\ tests are the same. This
is not a coincidence: if a UMPU test is unique (up to sets of measure zero)
and there exists a UMPI test with respect to some group of transformations,
then both coincide (up to sets of measure zero). For the IV model, however,
there are no uniformly most powerful tests. In perfect analogy to our
canonical model, there are two lines of research in the IV model. %
\citet{MoreiraMoreira19} seek optimal two-sided tests within a restricted
class of tests (the so-called SU\ tests) by fixing a long-run reduced-form
variance matrix, i.e., they consider the known and fixed case. In this
paper, we instead explore model symmetries by taking the reduced-form
variance to be known, but not fixed. As in the canonical model above, we
prefer not to take a stance on which thought experiment is more suitable. We
consider both approaches to be useful, leading to new insights in the IV
model.

If the error variance matrix in the instrumental variable regression is
considered known --but not fixed-- then the model satisfies some natural
symmetries. The main contribution of this paper is that we propose a test
that is invariant for the largest data transformation that leaves the model
and null hypothesis unchanged. The novel test, called the conditional
integrated likelihood (CIL) test, is invariant and the limit of WAP tests.
The weights are derived from relatively invariant measures on the parameter
space. The weights of the transformed parameters are then proportional to
the weights of the original parameters. The test statistic is the ratio of
the integrated likelihoods of the parameter space under the null and
alternative hypotheses. As a result, the invariance of the model combined
with the proportional effect of the transformation on the weights make the
CIL test invariant to the transformation, as required.

\subsection{Model-Preserving Transformations in the HAC-IV Model\label{Sec
Transformations}}

To understand model symmetries, it is convenient to transform the random
matrix $R$ into%
\begin{equation}
R_{0}=RB_{0}\text{, where }B_{0}=\left( 
\begin{array}{cc}
1 & 0 \\ 
-\beta _{0} & 1%
\end{array}%
\right) \text{,}  \label{R0}
\end{equation}%
so that the mean of the first column of $R_{0}=\left[ R_{1}:R_{2}\right] $
is zero under the null. The distribution of $R_{0}$ is%
\begin{equation}
R_{0}\sim N\left( \mu a_{\Delta }^{\prime },\Sigma _{0}\right) ,
\end{equation}%
where $a_{\Delta }^{\prime }=\left( \Delta ,1\right) $, $\Delta =\beta
-\beta _{0}$, and 
\begin{equation}
\Sigma _{0}=(B_{0}^{\prime }\otimes I_{k})\Sigma (B_{0}\otimes I_{k})=\left( 
\begin{array}{ll}
\Sigma _{11} & \Sigma _{12} \\ 
\Sigma _{21} & \Sigma _{22}%
\end{array}%
\right) .  \label{HAC}
\end{equation}%
We partition the inverse variance as%
\begin{eqnarray}
\Sigma _{0}^{-1} &=&\left[ 
\begin{array}{cc}
\Sigma ^{11} & \Sigma ^{12} \\ 
\Sigma ^{21} & \Sigma ^{22}%
\end{array}%
\right] \text{, where} \\
\Sigma ^{11} &=&\left( \Sigma _{11}-\Sigma _{12}\Sigma _{22}^{-1}\Sigma
_{21}\right) ^{-1}\text{, }\Sigma ^{22}=\left( \Sigma _{22}-\Sigma
_{21}\Sigma _{11}^{-1}\Sigma _{12}\right) ^{-1}\text{, and}  \notag \\
\Sigma ^{21} &=&\left( \Sigma ^{12}\right) ^{\prime }=-\Sigma ^{22}\Sigma
_{21}\Sigma _{11}^{-1}=-\Sigma _{22}^{-1}\Sigma _{21}\Sigma ^{11}\text{.} 
\notag
\end{eqnarray}

For $k=1$, the variance matrix $\Sigma $ trivially has a Kronecker
structure, as defined in Section \ref{Sec Kronecker}. Hence, AMS06 is
directly applicable. In particular, the Anderson-Rubin test is the UMPI test
in the just-identified model ($k=1$); see Comment 2 following Corollary 1 of
AMS06.\footnote{%
AMS06's optimality result for invariant tests when $k=1$ can be seen from
the perspective of unbiased tests. \citet{Moreira02,
Moreira09a} shows that the Anderson-Rubin test is uniformly most powerful
unbiased (UMPU). If there is a UMPI test, then the Anderson-Rubin test must
be the one; see Theorem 6.6.1 of \citet{LehmannRomano05}.}

For $k>1$, we recommend a novel WAP test. The weights are based on
invariance arguments. To show the model symmetries, we consider the affine
group of transformations $\left( A,G\right) \in R^{2k}\times R^{2k\times 2k}$
of $R_{0}$:%
\begin{equation}
A+G\cdot vec\left( R_{0}\right) \sim N\left( A+G\left( a_{\Delta }\otimes
\mu \right) ,G\Sigma _{0}G^{\prime }\right) \text{.}
\end{equation}%
If we consider $\Sigma _{0}$ to be fixed, we have to impose restrictions on $%
G$ and/or $\Sigma _{0}$ for the transformation to preserve the model, so
that 
\begin{equation}
G\Sigma _{0}G^{\prime }=\Sigma _{0}.
\end{equation}%
If the variance matrix $\Sigma _{0}$ is known but not fixed, it changes with
the transformation. If $\Sigma _{0}$ is a known variance matrix, so is $%
G\Sigma _{0}G^{\prime }$.

Partitioning $A$ into $k$-dimensional vectors and $G$ into $k\times k$
matrices: 
\begin{equation}
A=\left[ 
\begin{array}{c}
A_{1} \\ 
A_{2}%
\end{array}%
\right] \text{ and }G=\left[ 
\begin{array}{cc}
G_{11} & G_{21} \\ 
G_{21} & G_{22}%
\end{array}%
\right] ,
\end{equation}%
we find that the expectation of the transformed $R_{0}$ becomes 
\begin{equation}
\mathbb{E}[A+G\cdot vec\left( R_{0}\right) ]=A+\left[ 
\begin{array}{cc}
G_{11} & G_{12} \\ 
G_{21} & G_{22}%
\end{array}%
\right] \cdot \left[ 
\begin{array}{c}
\Delta \mu \\ 
\mu%
\end{array}%
\right] =\left[ 
\begin{array}{c}
A_{1}+\left( G_{11}\Delta +G_{12}\right) \mu \\ 
A_{2}+\left( G_{21}\Delta +G_{22}\right) \mu%
\end{array}%
\right] .
\end{equation}%
To preserve the null hypothesis $H_{0}:\Delta =0$, the first sub-vector of
the mean has to be zero for all values of $\mu $. This forces $A_{1}=0$ and $%
G_{12}=0$. In the original model, the mean of $R_{1}$ is proportional to the
mean of $R_{2}$. To preserve the model, the two subvectors of the
transformed mean must be proportional to each other, which forces $A_{2}=0$
and 
\begin{equation}
G_{11}\Delta \cdot \mu \varpropto \left( G_{21}\Delta +G_{22}\right) \mu
\end{equation}%
for all $\mu $. This implies that $G_{11}=g_{11}\cdot
g_{1},G_{21}=g_{21}\cdot g_{1},G_{22}=g_{22}\cdot g_{1}$ with $%
g_{11},g_{21},g_{22}$ being constants of proportionality, and $g_{1}$ a $%
k\times k$ matrix. As a result,%
\begin{equation}
G=\left[ 
\begin{array}{cc}
g_{11}\cdot g_{1} & 0 \\ 
g_{21}\cdot g_{1} & g_{22}\cdot g_{1}%
\end{array}%
\right] =g_{2}\otimes g_{1},
\end{equation}%
where $g_{2}$ is a $2\times 2$ lower triangular matrix. Therefore, $g=\left(
g_{1},g_{2}\right) $ transforms the data to 
\begin{equation}
g\circ \left( R_{0},\Sigma _{0}\right) =\left( g_{1}R_{0}g_{2}^{\prime
},\left( g_{2}\otimes g_{1}\right) \Sigma _{0}\left( g_{2}^{\prime }\otimes
g_{1}^{\prime }\right) \right) .
\label{(action sample-space non-Kronecker 3)}
\end{equation}%
We use the transpose of $g_{2}$ so that the associated transformation is a 
\textit{left action}. This transformation leaves the model unchanged: it
preserves the null and the proportionality of the subvectors of the mean of $%
R_{0}$. Specifically,%
\begin{equation}
g\circ \left( \Delta ,\mu ,\Sigma _{0}\right) =\left( \frac{\Delta g_{11}}{%
\Delta g_{21}+g_{22}},g_{1}\mu \left( \Delta g_{21}+g_{22}\right) ,\left(
g_{2}\otimes g_{1}\right) \Sigma _{0}\left( g_{2}^{\prime }\otimes
g_{1}^{\prime }\right) \right) .
\label{(left action
parameter-space non-Kronecker)}
\end{equation}%
If the matrix $\Sigma _{0}$ is invertible, as we assume here, then $g_{1}$
and $g_{2}$ must be non-singular.\footnote{%
Therefore, $g_{11}$ and $g_{22}$ are non-zero elements.} This means that $%
g_{1}\in \mathcal{G}_{L}\left( k\right) $ and $g_{2}\in \mathcal{G}%
_{T}\left( 2\right) $, respectively, the groups of invertible $k\times k$
matrices and of invertible lower-triangular $2\times 2$ matrices (with
matrix multiplication as the group operator).

If the original data are supportive of the null hypothesis, then the
transformed data should be equally supportive of this hypothesis. The test
should be the same whether it is computed from the original or from the
transformed data, i.e. the test should be invariant to the transformation $g$%
. As we show later, the AR, LM, CQLR, and CLR tests are invariant to the
group of transformations $g$ presented above. Without further restrictions
on the weight $m\left( T\right) $, the CLC test may be sensitive to this
transformation. Likewise, the weighted-average-power SU test proposed by %
\citet{MoreiraMoreira19} may also change with data transformations. As a
result, the CLC and SU tests may have power that changes with the data
transformation. Next, we propose a conditional integrated weighted
likelihood test that is invariant to $g$. This test does not have the
undesirably low power of WAP tests based on generic weights documented by %
\citet{MoreiraMoreira19}.

\subsection{The CIL\ Test\label{Sec CIL Test}}

Consider an integrated likelihood ($IL$)\ statistic which is the ratio of
two terms. The numerator is the integrated likelihood over $\mu $ with
respect to the Lebesgue measure and over $\Delta $ with respect to $%
\left\vert \Delta \right\vert ^{k-2}d\Delta $. The denominator is the
density of the pivotal statistic $S$ under the null hypothesis. In Appendix
A, we show the $IL$ statistic is%
\begin{eqnarray}
IL &=&\int_{-\infty }^{\infty }e^{\frac{1}{2}\left[ vec\left( R_{0}\right)
^{\prime }\Sigma _{0}^{-1/2}N_{\Sigma _{0}^{-1/2}(a_{\Delta }\otimes
I_{k})}\Sigma _{0}^{-1/2}vec\left( R_{0}\right) -T^{\prime }T\right] }
\label{(IL stat)} \\
&&\times \left\vert \left( a_{\Delta }^{\prime }\otimes I_{k}\right) \Sigma
_{0}^{-1}\left( a_{\Delta }\otimes I_{k}\right) \right\vert
^{-1/2}\left\vert \Delta \right\vert ^{k-2}d\Delta ,  \notag
\end{eqnarray}%
up to a multiplication by $\left\vert \Sigma _{0}\right\vert ^{1/2}$. We
also prove this integral is finite in the over-identified case $k\geq 2$.
The conditional (on $T$) integrated likelihood (CIL) test based on (\ref{(IL
stat)}) is invariant to the transformation $g$ because both the $IL$
statistic and its conditional quantile have the same proportionality
multiplier $\chi \left( g\right) $ with respect to $g$. Furthermore, the
CIL\ test is the limit of a sequence of WAP tests. We relegate the theory
and proofs to Section \ref{Sec Invariant Tests}. Here, we focus on the
implementation of the CIL\ test.

The integral defined in (\ref{(IL stat)}) is improper, which can create
computational difficulties. We circumvent this problem by changing
variables, so that the integral is proper. This is convenient for the
numerical integration that we use to compute the $IL$ statistic. First, we
standardize the vector $a_{\Delta }$ to have norm one: 
\begin{equation}
\overline{a}_{\Delta }=\frac{a_{\Delta }}{\left( 1+\Delta ^{2}\right) ^{1/2}}%
=\left( \frac{\Delta }{\left( 1+\Delta ^{2}\right) ^{1/2}},\frac{1}{\left(
1+\Delta ^{2}\right) ^{1/2}}\right) .
\end{equation}%
We note that $N_{\Sigma _{0}^{-1/2}(a_{\Delta }\otimes I_{k})}=N_{\Sigma
_{0}^{-1/2}(\overline{a}_{\Delta }\otimes I_{k})}$ and also that%
\begin{equation}
\left\vert \left( a_{\Delta }^{\prime }\otimes I_{k}\right) \Sigma
_{0}^{-1}\left( a_{\Delta }\otimes I_{k}\right) \right\vert ^{-1/2}=\left(
1+\Delta ^{2}\right) ^{-k/2}\left\vert \left( \overline{a}_{\Delta }^{\prime
}\otimes I_{k}\right) \Sigma _{0}^{-1}\left( \overline{a}_{\Delta }\otimes
I_{k}\right) \right\vert ^{-1/2}.
\end{equation}%
Therefore,%
\begin{eqnarray}
IL &=&\int_{-\infty }^{\infty }e^{\frac{1}{2}\left[ vec\left( R_{0}\right)
^{\prime }\Sigma _{0}^{-1/2}N_{\Sigma _{0}^{-1/2}(\overline{a}_{\Delta
}\otimes I_{k})}\Sigma _{0}^{-1/2}vec\left( R_{0}\right) -T^{\prime }T\right]
} \\
&&\times \left\vert \left( \overline{a}_{\Delta }^{\prime }\otimes
I_{k}\right) \Sigma _{0}^{-1}\left( \overline{a}_{\Delta }\otimes
I_{k}\right) \right\vert ^{-1/2}\left\vert \frac{\Delta }{\left( 1+\Delta
^{2}\right) ^{1/2}}\right\vert ^{k}\frac{1}{\Delta ^{2}}d\Delta .  \notag
\end{eqnarray}

By changing variables following 
\begin{equation}
l_{\eta }\equiv \left( \sin \eta ,\cos \eta \right) ^{\prime }=\overline{a}%
_{\Delta },
\end{equation}%
the $IL$ statistic becomes 
\begin{eqnarray}
IL &=&\int_{-pi/2}^{pi/2}e^{\frac{1}{2}\left[ vec\left( R_{0}\right)
^{\prime }\Sigma _{0}^{-1/2}N_{\Sigma _{0}^{-1/2}(l_{\eta }\otimes
I_{k})}\Sigma _{0}^{-1/2}vec\left( R_{0}\right) -T^{\prime }T\right] } \\
&&\times \left\vert \left( l_{\eta }^{\prime }\otimes I_{k}\right) \Sigma
_{0}^{-1}\left( l_{\eta }\otimes I_{k}\right) \right\vert ^{-1/2}\left\vert
\sin \eta \right\vert ^{k}\frac{\left( 1+\tan ^{2}\eta \right) }{\tan
^{2}\eta }d\eta ,  \notag
\end{eqnarray}%
where $pi=3.14159...$. Because%
\begin{equation}
\left\vert \sin \eta \right\vert ^{k}\frac{\left( 1+\tan ^{2}\eta \right) }{%
\tan ^{2}\eta }=\left\vert \sin \eta \right\vert ^{k-2}\frac{\sin \eta ^{2}}{%
\cos ^{2}\eta \cdot \tan ^{2}\eta }=\left\vert \sin \eta \right\vert ^{k-2},
\end{equation}%
the $IL$ statistic simplifies to%
\begin{eqnarray}
IL &=&\int_{-pi/2}^{pi/2}e^{\frac{1}{2}\left[ vec\left( R_{0}\right)
^{\prime }\Sigma _{0}^{-1/2}N_{\Sigma _{0}^{-1/2}(l_{\eta }\otimes
I_{k})}\Sigma _{0}^{-1/2}vec\left( R_{0}\right) -T^{\prime }T\right] }
\label{IL-stat} \\
&&\times \left\vert \left( l_{\eta }^{\prime }\otimes I_{k}\right) \Sigma
_{0}^{-1}\left( l_{\eta }\otimes I_{k}\right) \right\vert ^{-1/2}\left\vert
\sin \eta \right\vert ^{k-2}d\eta ,  \notag
\end{eqnarray}%
which is easier to compute.

The $IL$ statistic can be compared to the statistic that integrates the
likelihood with respect to the Lebesgue measure $d\mu \times d\Delta $
without the weights $|\Delta |^{k-2}$:%
\begin{eqnarray}
IL_{0} &=&\int_{-pi/2}^{pi/2}e^{\frac{1}{2}\left[ vec\left( R_{0}\right)
^{\prime }\Sigma _{0}^{-1/2}N_{\Sigma _{0}^{-1/2}(l_{\eta }\otimes
I_{k})}\Sigma _{0}^{-1/2}vec\left( R_{0}\right) -T^{\prime }T\right] }
\label{(IL polar)} \\
&&\times \left\vert \left( l_{\eta }^{\prime }\otimes I_{k}\right) \Sigma
_{0}^{-1}\left( l_{\eta }\otimes I_{k}\right) \right\vert ^{-1/2}\left\vert
\cos \eta \right\vert ^{k-2}d\eta .  \notag
\end{eqnarray}%
(The derivation of $IL_{0}$ is analogous to the $IL$ statistic, as shown in
Appendix A.) Numerically, the computation of $IL$ or $IL_{0}$ is equally
difficult. Without the weights $|\Delta |^{k-2}$, the statistic $IL_{0}$
does not yield an invariant test when $k>2$. Hence, the test suffers the
power problems documented by \citet{MoreiraMoreira19}.

The representation of $IL$ in terms of $R_{0}$ and $\Sigma _{0}$ is
convenient to prove that the CIL\ test is invariant to the transformation (%
\ref{(action sample-space non-Kronecker 3)}). However, the approach can be
unnecessarily challenging when testing for different levels of $\beta _{0}$.
This pitfall can be important to derive confidence regions, which consist of
all values of $\beta _{0}$ which are not rejected by the CIL test. For
numerical stability, we instead recommend representing $IL$ in terms of the
original data, $R$ and $\Sigma $.

Algebraic manipulations show that%
\begin{eqnarray}
IL &=&\int_{-\infty }^{\infty }e^{\frac{1}{2}\left[ vec\left( R\right)
^{\prime }\Sigma ^{-1/2}N_{\Sigma ^{-1/2}(a\otimes I_{k})}\Sigma
^{-1/2}vec\left( R\right) -T^{\prime }T\right] } \\
&&\times \left\vert \left( a^{\prime }\otimes I_{k}\right) \Sigma
^{-1}\left( a\otimes I_{k}\right) \right\vert ^{-1/2}\left\vert \beta -\beta
_{0}\right\vert ^{k-2}d\beta .  \notag
\end{eqnarray}%
By changing variables 
\begin{equation}
l_{\theta }\equiv \left( \sin \theta ,\cos \theta \right) ^{\prime }=\frac{a%
}{\left\Vert a\right\Vert }\equiv \overline{a},
\end{equation}%
and following steps analogous to the derivation of (\ref{IL-stat}), we show
in Appendix A that%
\begin{eqnarray}
IL &=&\left( 1+\beta _{0}^{2}\right) ^{\left( k-2\right)
/2}\int_{-pi/2}^{pi/2}e^{\frac{1}{2}\left[ vec\left( R\right) ^{\prime
}\Sigma ^{-1/2}N_{\Sigma ^{-1/2}(l_{\theta }\otimes I_{k})}\Sigma
^{-1/2}vec\left( R\right) -T^{\prime }T\right] } \\
&&\times \left\vert \left( l_{\theta }^{\prime }\otimes I_{k}\right) \Sigma
^{-1}\left( l_{\theta }\otimes I_{k}\right) \right\vert ^{-1/2}\left\vert 
\frac{1}{\sqrt{1+\beta _{0}^{2}}}\sin \theta -\frac{\beta _{0}}{\sqrt{%
1+\beta _{0}^{2}}}\cos \theta \right\vert ^{k-2}d\theta .  \notag
\end{eqnarray}%
The factor $\left( 1+\beta _{0}^{2}\right) ^{\left( k-2\right) /2}$ can be
ignored in the computation of the CIL\ test, as it is directly absorbed by
the critical value function. Hence, we suggest implementing the conditional
test based on the $\left( 1+\beta _{0}^{2}\right) ^{-\left( k-2\right) /2}IL$
statistic.

There are also connections between the $IL$ statistic and the $LR$
statistic. The $LR$ statistic maximizes, with respect to $\Delta $,%
\begin{equation}
\left[ vec\left( R_{0}\right) ^{\prime }\Sigma _{0}^{-1/2}N_{\Sigma
_{0}^{-1/2}(a_{\Delta }\otimes I_{k})}\Sigma _{0}^{-1/2}vec\left(
R_{0}\right) -T^{\prime }T\right] ,
\end{equation}%
which is the term inside the brackets of (\ref{(IL stat)}). The $IL$
statistic integrates the exponential of this term after two corrections. The
first correction, $\left\vert \left( l_{\eta }^{\prime }\otimes I_{k}\right)
\Sigma _{0}^{-1}\left( l_{\eta }\otimes I_{k}\right) \right\vert ^{-1/2}$,
arises from integration with respect to the Lebesgue measure $d\mu $. The
second correction $\left\vert \sin \eta \right\vert ^{k-2}$ ensures that the
test is two-sided and invariant, so that we avoid the one-sided power
behavior in parts of the parameter space. In the next section, we show some
advantages of the CIL\ test over the AR and CLR tests.

\section{Numerical Simulations \label{Sec Numerical}}

Here, we provide numerical simulations for the AR, CLR, CIL, and CIL$_{0}$\
tests. All results reported here are for $k=5$ and only one level of
instrument strength based on 1,000 Monte Carlo replications for power and
1,000 simulations to approximate the tests' critical value function. In the
supplement, we provide two levels of identification strength and consider $%
k=2,5,10$. For reasons explained below, a reliable implementation for the
CLR test is computationally intensive. Because of this, the supplemental
power plots are limited to only 200 replications and 200 simulations for the
conditional quantile.

We first illustrate numerical problems with likelihood optimization and
integration. Some of these difficulties arise even in the simple case in
which errors are homoskedastic. We focus on tests with significance level
5\% for testing $\beta _{0}=0$. We set the parameter $\mu =\left( \lambda
^{1/2}/\sqrt{k}\right) 1_{k}$ for $k=5$ and set $\lambda /k=2$, where $1_{k}$
is a $k$-dimensional vector of ones and $\lambda $ is a measure of the IVs'
strength. The variance of structural-form errors is one and their
correlation is $\rho =-0.9,0.9$. We present plots for the power envelope and
power functions against various alternative values of $\beta $. We plot
power as a function of the rescaled alternative $\beta \lambda ^{1/2}$,
which reflects the difficulty of making inference on $\beta $ for different
instruments' strength. 
\begin{figure}[th]
\caption{Power Curves for Homoskedastic Errors and $k=5$ \newline
{(Likelihood Optimization)}}
\label{fig: CLR homo k=5}%
\begin{subfigure}{.5\textwidth} \centering
\includegraphics[width=.9\linewidth]{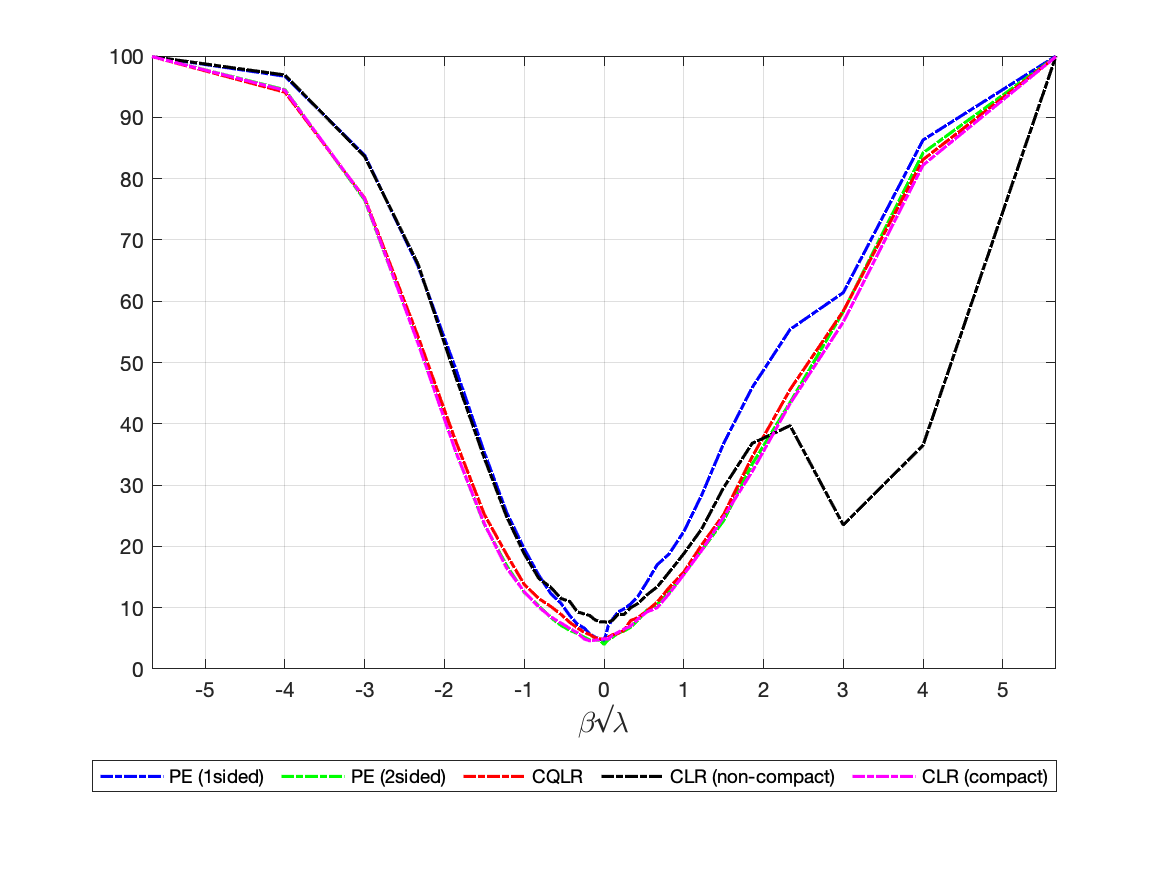}
\subcaption{$\rho=0.9$}
\end{subfigure}
\begin{subfigure}{.5\textwidth}
\centering
\includegraphics[width=.9\linewidth]{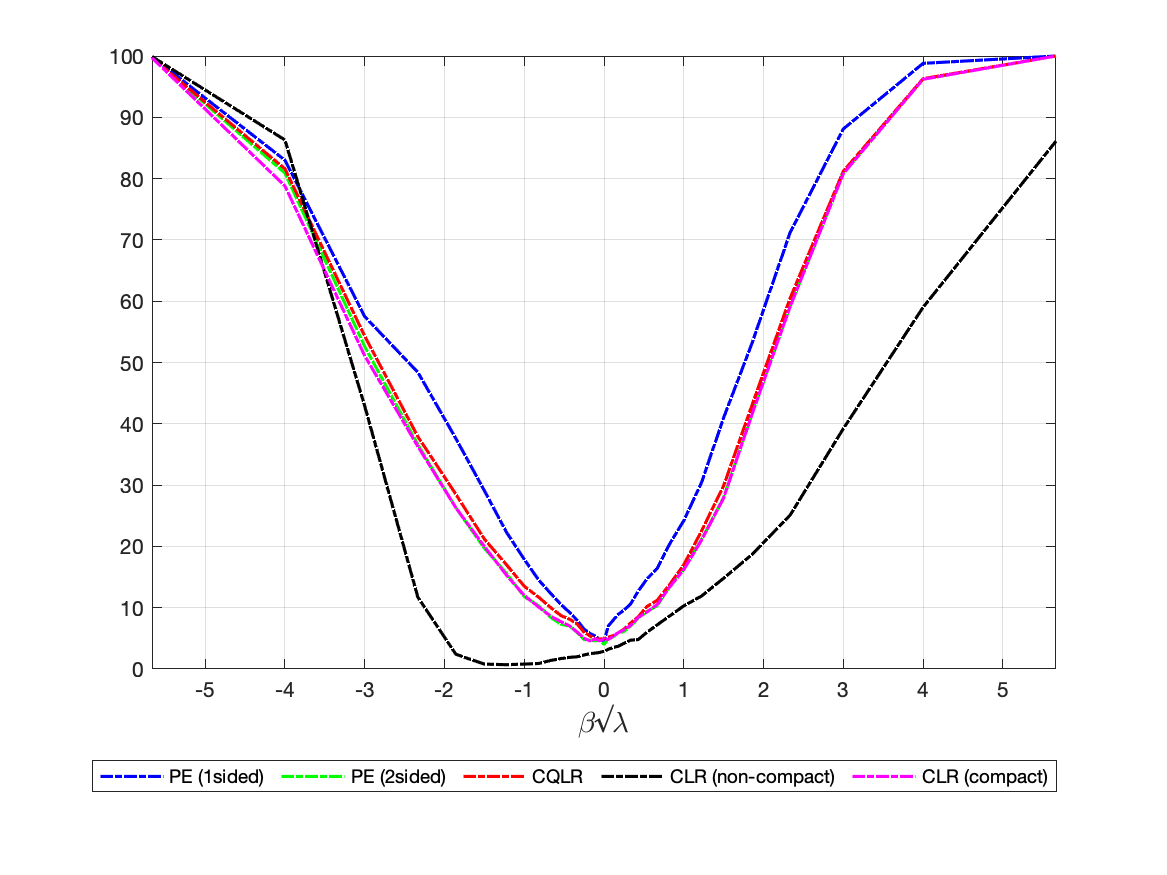}
\subcaption{$\rho=-0.9$}
\end{subfigure}
\end{figure}

Figure \ref{fig: CLR homo k=5} presents the one-sided and two-sided power
envelope for invariant similar tests. These power envelopes are derived
analytically by \citet{MillsMoreiraVilela14} and %
\citet{AndrewsMoreiraStock06}, respectively; see earlier theory by %
\citet{AndrewsMoreiraStock04}. This early work also shows these power bounds
are valid for all invariant tests which have correct size. We also plot
power curves for the CQLR test as well as two numerical optimization
strategies to obtain the CLR test. The first randomly draws the initial
point for the search optimization algorithm in (\ref{(LR stat)}) for $\beta $%
. Here, we consider the uniform distribution over $\left[ -1000,1000\right] $%
. The second one relies on the fact that we can maximize the likelihood over
a compact set, without loss of generality. We can write the likelihood ratio
statistic as%
\begin{equation}
LR=\max_{\theta }vec(R)^{\prime }\Sigma ^{-1/2}N_{\Sigma ^{-1/2}(l_{\theta
}\otimes I_{k})}\Sigma ^{-1/2}vec(R)-T^{\prime }T,
\end{equation}%
where the maximization is over the compact set $\left[ -pi/2,pi/2\right] $.
The initial point is drawn from a uniform distribution over that same set.
Recall that the CQLR and CLR tests are theoretically identical when errors
are homoskedastic. \emph{Any} power difference between the CQLR test and
these numerical implementations for the CLR test arises from failures in the
likelihood optimization.

The power upper bounds are useful to understand the difficulty in the
likelihood optimization behind the CLR test. Both CQLR and CLR tests based
on optimization over the compact set for $\theta $ perform alike. These two
tests have power very close to the two-sided power envelope. The CLR test
based on a draw-and-search for the optimal $\beta $ fails remarkably. In the
first graph, this implementation has power above the two-sided power
envelope and close to the one-sided power bound for parts of the parameter
space. Furthermore, this implementation must fail to deliver a test with
correct size. Indeed, the implementation for the CLR test over the whole
real line has size close to 10\% instead of the correct 5\% level. To make
matters even worse, the second plot in Figure \ref{fig: CLR homo k=5} shows
bad behavior associated with the sample implementation of the CLR test. The
power can even be close to zero for parts of the parameter space.

Of course, one could simply use the CQLR test for the homoskedastic case.
The lesson learned here is that likelihood optimization does matter for the
power performance of the CLR test in general. In more complex designs (i.e.,
non-Kronecker error variance), drawing a unique initial point is far from
sufficient. Our experience is that likelihood maximization for the
implementation of the CLR test can be very slow and unreliable. This is
particularly true when several initial points are required, as happens in
some designs below. 
\begin{figure}[th]
\caption{Power Curves for Homoskedastic Errors and $k=5$ \newline
{(Integrated Likelihood)}}
\label{fig: CIL homo k=5}%
\begin{subfigure}{.5\textwidth} \centering
\includegraphics[width=.9\linewidth]{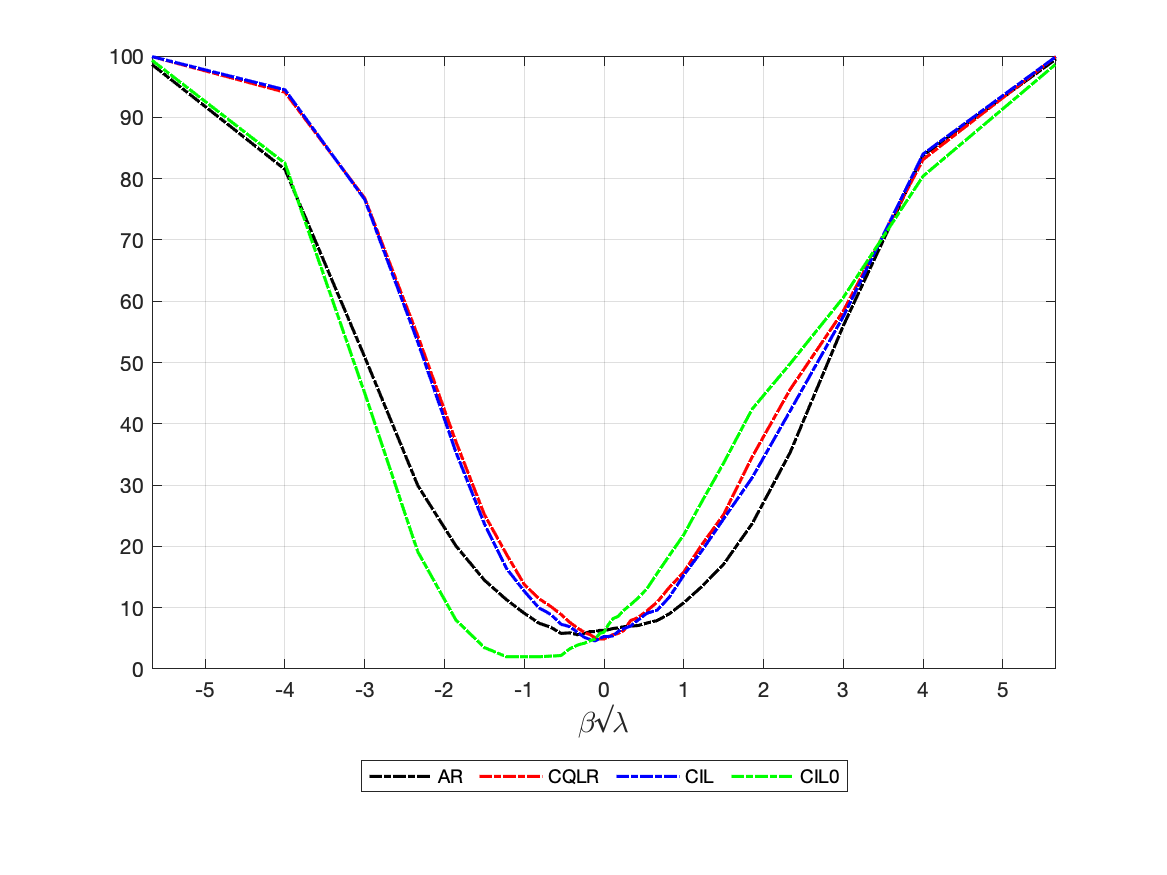}
\subcaption{$\rho=0.9$}
\end{subfigure}
\begin{subfigure}{.5\textwidth}
\centering
\includegraphics[width=.9\linewidth]{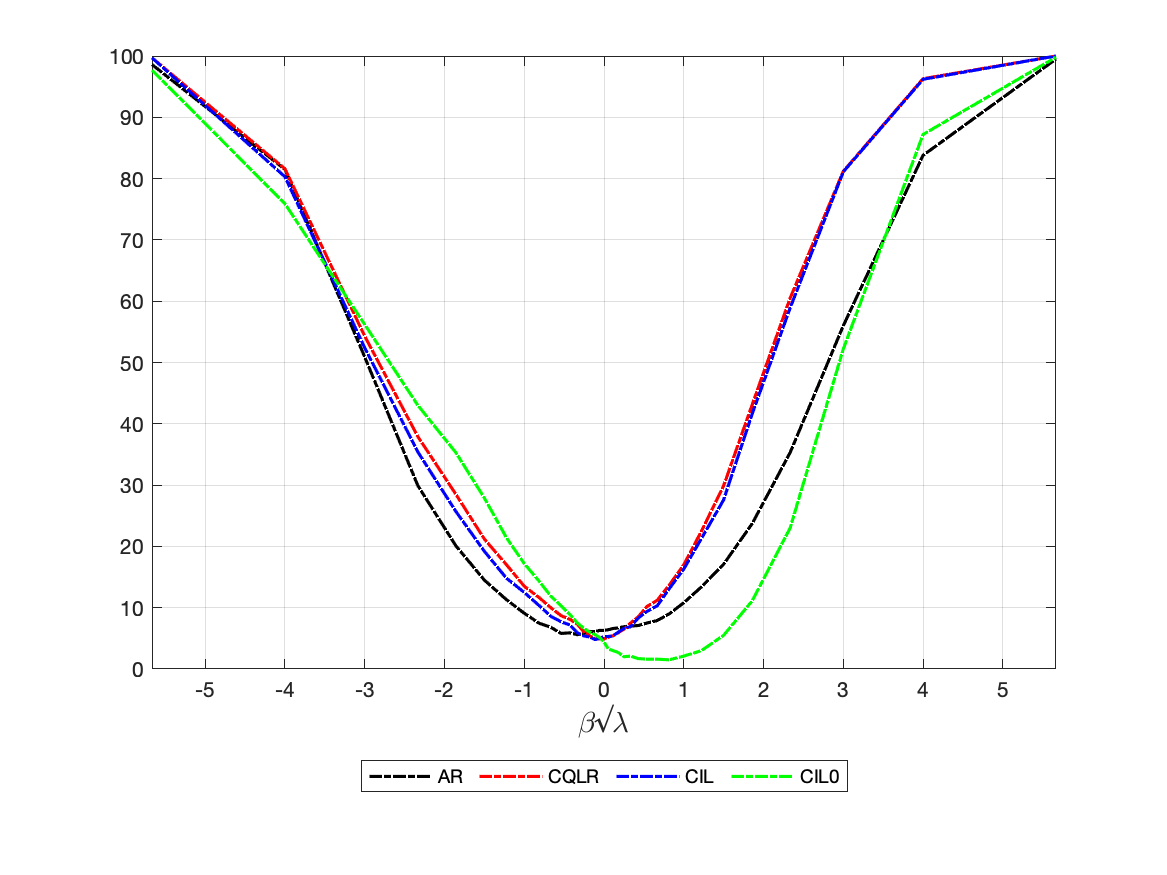}
\subcaption{$\rho=-0.9$}
\end{subfigure}
\end{figure}

Figure \ref{fig: CIL homo k=5} presents power for the AR, CQLR, CIL, and CIL$%
_{0}$ tests. The CQLR and CIL\ tests have comparable power and outperform
the AR test. These plots are a reassurance that the CIL\ test performs well
in scenarios more favorable to CQLR. The CIL$_{0}$ test has behavior quite
different from the CIL\ test. The dissimilar behavior of the CIL\ and CIL$%
_{0}$ tests illustrates that tests based on likelihood integration are
sensitive to weight choices. While the CIL\ and CIL$_{0}$ tests perform
comparably when $\lambda $ increases, they have distinct properties when IVs
are weak. For one side of the alternative, the power of CIL$_{0}$ is smaller
than that of the CQLR and CIL\ tests. For the other side of the alternative,
it actually has larger power. Therefore, the CIL$_{0}$ behaves as a
one-sided test. The CIL$_{0}$ test being biased means the null rejection
probability is smaller for some alternatives than under the null. This
undesirable feature of the CIL$_{0}$ test is not shared by the CQLR and CIL
test. These two tests do not suffer the same power deficiencies as the CIL$%
_{0}$ test. They behave as two-sided tests by construction, and have power
close to the power upper bound.

We now move to the more complex case in which errors can be heteroskedastic,
autocorrelated, and/or clustered (HAC). We replicate four designs: the
near-singular (NS), a variation thereof (NS with perturbation), and growing
alternative (GA) designs of \citet{MoreiraRidder20}, and the non-Kronecker
(NK) design of \citet{MoreiraMoreira19}. While these simulations are not
exhaustive for all parameter combinations, none of these designs is chosen
to favor the CIL\ test over the CLR test. The main goal of these designs is
only to show that there exist invariant tests which depend on the statistic $%
S$ beyond $AR$ and $LM$.

\begin{table}[tbph]
\caption{Likelihood Optimization and Initial Values (percentage)}
\label{tab:Initial Values Percentage}\centering
\begin{tabular}{l|cccc}
\hline\hline
& (1,No) & (0,Yes) & (51,No) & (50, Yes) \\ \hline
(1,No) & - & 64.7 & 86.7 & 87.7 \\ 
(0,Yes) & 26.0 & - & 47.1 & 47.2 \\ 
(51,No) & 0.1 & 3.5 & - & 4.2 \\ 
(50,Yes) & 0.1 & - & 0.4 & - \\ \hline
\end{tabular}%
\end{table}

To conserve space, we focus here only on simulations based on the NS design
for $k=5$. We set $\mu =\lambda ^{1/2}e_{1}$, with $\lambda /k=2$. For the
variance matrix, we define $J_{k}$ to be the $k\times k$ matrix with the
anti-diagonal elements equal to one and the other components zero. We have $%
J_{k}^{2}=I_{k}$. The $k\times k$ submatrices of $\Sigma _{0}$ are%
\begin{equation}
\Sigma _{11}=c_{11}I_{k}\text{, }\Sigma _{12}=c_{12}J_{k}\text{, and }\Sigma
_{22}=c_{22}I_{k}\text{,}
\end{equation}%
where $c_{11}$, $c_{12}$, and $c_{22}$ are tuning parameters. The values for
the NS design are $c_{11}=1$, $c_{12}=100$, and $%
c_{22}=c_{12}^{2}+c_{12}^{-3}$. In this design, the power of both LM and
CQLR tests is essentially equal to size. The full set of results for $%
k=2,5,10$ and $\lambda /k=2,10$ for all four designs, as well as
descriptions of the NS with perturbation, GA, and NK designs, are presented
in the supplement. 
\begin{table}[tbph]
\caption{Likelihood Optimization and Initial Values (factor)}
\label{tab:Initial Values Factor}\centering
\begin{tabular}{l|cccc}
\hline\hline
& (1,No) & (0,Yes) & (51,No) & (50, Yes) \\ \hline
(1,No) & - & 3,270.1 & 5,446.5 & 5,410.0 \\ 
(0,Yes) & 964.9 & - & 9,910.7 & 9,890.0 \\ 
(51,No) & 1,114.7 & 452.6 & - & 379.1 \\ 
(50,Yes) & 0.2 & - & 7.7 & - \\ \hline
\end{tabular}%
\end{table}

When the variance matrix has a Kronecker product form, the $LR$ statistic
has a closed-form solution, and the CLR test reduces to the CQLR test. This
sidesteps the daunting task of numerically optimizing the likelihood. In the
special case with homoskedastic errors, choosing only one initial point
after compactifying the search set is enough for our purposes.
Unfortunately, this conclusion is not valid for more complex variance
matrices. Table \ref{tab:Initial Values Percentage} assesses improvements
for the likelihood optimization under the null hypothesis. The values inside
the parentheses indicate the number of random initial points for $\theta $
and whether $\beta _{0}$ is included or not, respectively. We compute the LR
statistic over 1,000 simulations for each case. We then report the
proportion of times in which one setup outperforms another setup (relative
improvement by an error margin of at least 0.1\%).

Each row in Table \ref{tab:Initial Values Percentage} corresponds to a
choice of the number of starting values and whether $\beta _{0}$ is among
the starting values, as specified in the row header. The entries in a row
report the fraction of repetitions in which the initial values selection and
the inclusion of $\beta _{0}$, as specified in the column header, give a
higher maximum likelihood value. For example, if we choose $\beta _{0}$
instead of only one random point as the initial value, we see improvements
in the likelihood optimization 64.7\% of the time. Conversely, the
likelihood optimization performs better 26.0\% of the time if we choose a
random point instead of $\beta _{0}$. For both of these scenarios,
improvements are gained by adding about 50 random initial values. This can
be seen in the upper-right $2\times 2$ block in Table \ref{tab:Initial
Values Percentage}, where the improvements range from 47.1\% to 87.7\%. On
the other hand, the improvements are negligible from starting with 50 random
points and $\beta _{0}$ as initial values --even when we include 51 other
random points. What is perhaps interesting is the improvement of 3.5\% from
adding $\beta _{0}$ as an initial value in addition to the 50 random points.
These two findings suggest that running optimization algorithms after
including 50 random points and $\beta _{0}$ as initial values should suffice
for our purposes. More worrisome, for smaller values of $\lambda $ or other
combinations of $\mu $ and the variance matrix, we may need to include even
more initial points. This may happen, for example, if the likelihood can be
flat for parts of the parameter space.

Table \ref{tab:Initial Values Factor} presents the average percentage
improvement (for the observations in which the error margin is at least
0.1\%). Even when we include 51 random initial points, meaningful
improvements can be gained by including the unknown parameter $\beta $.
These gains are on the order of 379.1\% for 4.2\% of the replications when
we include another 50 random initial points and $\beta $ itself. On the
other hand, when we include 51 other random initial points beyond $\beta $
and 50 points, the average improvement is on the order of 7.7\% for only
0.4\% of the repetitions.

All simulations are for the null hypothesis. For the alternative, it is
natural to use 50 random points, the null $\beta _{0}$, and the alternative $%
\beta $ as initial points. Of course, the parameter $\beta $ is unknown.
However, we want to minimize the numerical issues associated with the CLR
test, in case better optimization methods are found in the future. The table
shows that the solution of using $\beta $ in addition to 50 random points
works well to compute the LR statistic. A more complex problem happens when
we find the approximation for the critical value function. Recall that this
function is the conditional quantile under the null hypothesis. This
quantile is found by generating $S$ from a standard multivariate normal
distribution. That means the model is misspecified when $T$ is not generated
under the null. One possibility is to use the pseudo-parameter which
minimizes the Kullback-Leibler divergence criterion. This strategy follows
from the fact that the maximum likelihood (ML) estimator converges to this
pseudo-parameter under strong instruments. This route seems complicated and
unnecessary for our purposes. Excluding this parameter, we get smaller
values for the test statistic --not larger. Hence, the 95\% quantile used
for the critical value function tends to underestimate the true conditional
quantile. The bottom line is that by including $\beta $, the $LR$ statistic
is optimized properly while the conditional quantile can be smaller than it
should be. This means that, if anything, we may be overestimating the power
of the CLR test.

Finally, we evaluate improvements over other numbers of random initial
points for $\theta $. For example, unreported simulations show gains of
about 5.3\% obtained from adding 1 initial random point beyond 20 random
initial points. The choice of 50 seems the most sensible, in terms of
reliability and computational speed. Even then, the computation time for CLR
is about 35 times slower than that of the CIL\ test, on average (with the
range between 4 to 100 times slower). For the aforementioned reasons, we
include $\beta _{0}$ and $\beta $ as initial points as well. At least for
the designs considered here, unreported power comparisons for different
choices of initial values indicate that including 50 random points, $\beta
_{0}$, and $\beta $ offers\ stable and reliable power curves. There is, of
course, no guarantee that this searching scheme would be sufficient for
other designs.

\begin{figure}[th]
\caption{Power Curves for HAC Errors with $k=5$ and $\protect\lambda /k=2$}
\label{fig: HAC k=5}%
\begin{subfigure}{.5\textwidth} \centering
\includegraphics[width=.9\linewidth]{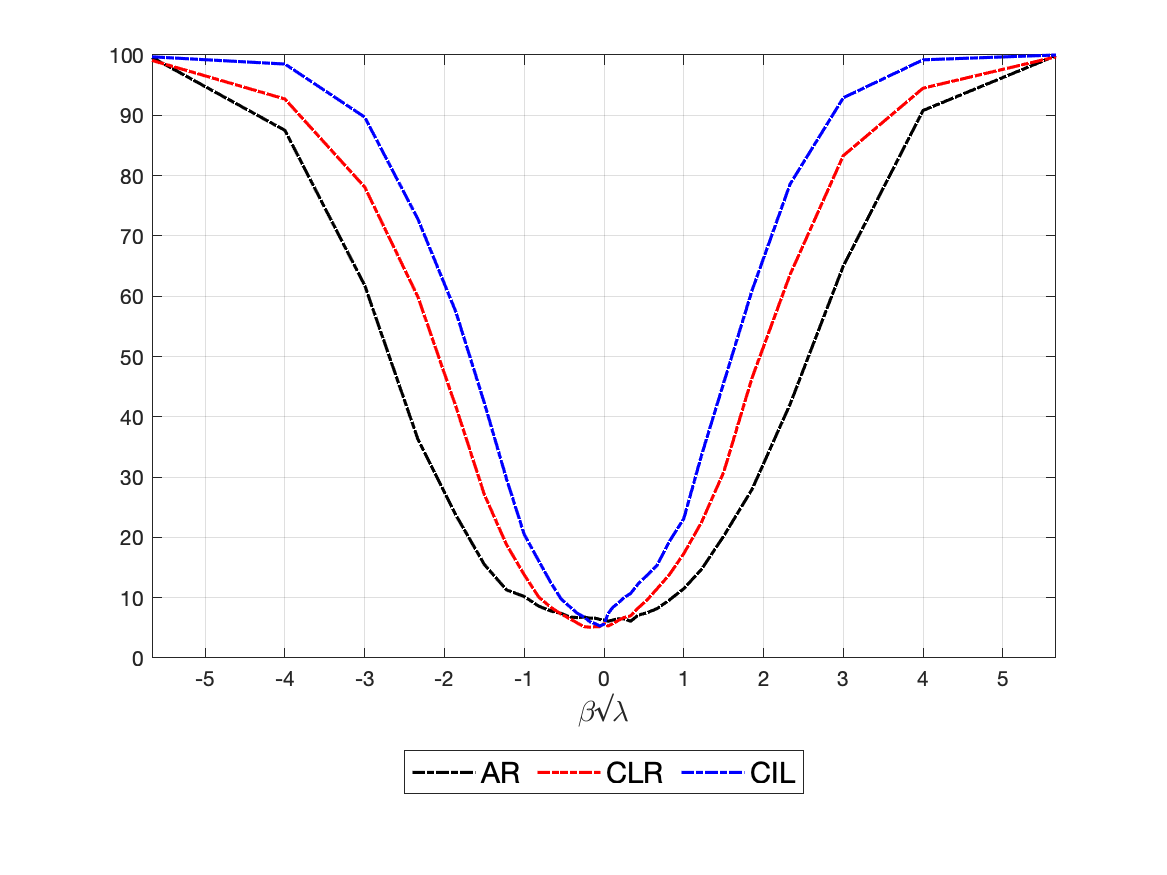}
\subcaption{NS Design}
\end{subfigure}
\begin{subfigure}{.5\textwidth}
\centering
\includegraphics[width=.9\linewidth]{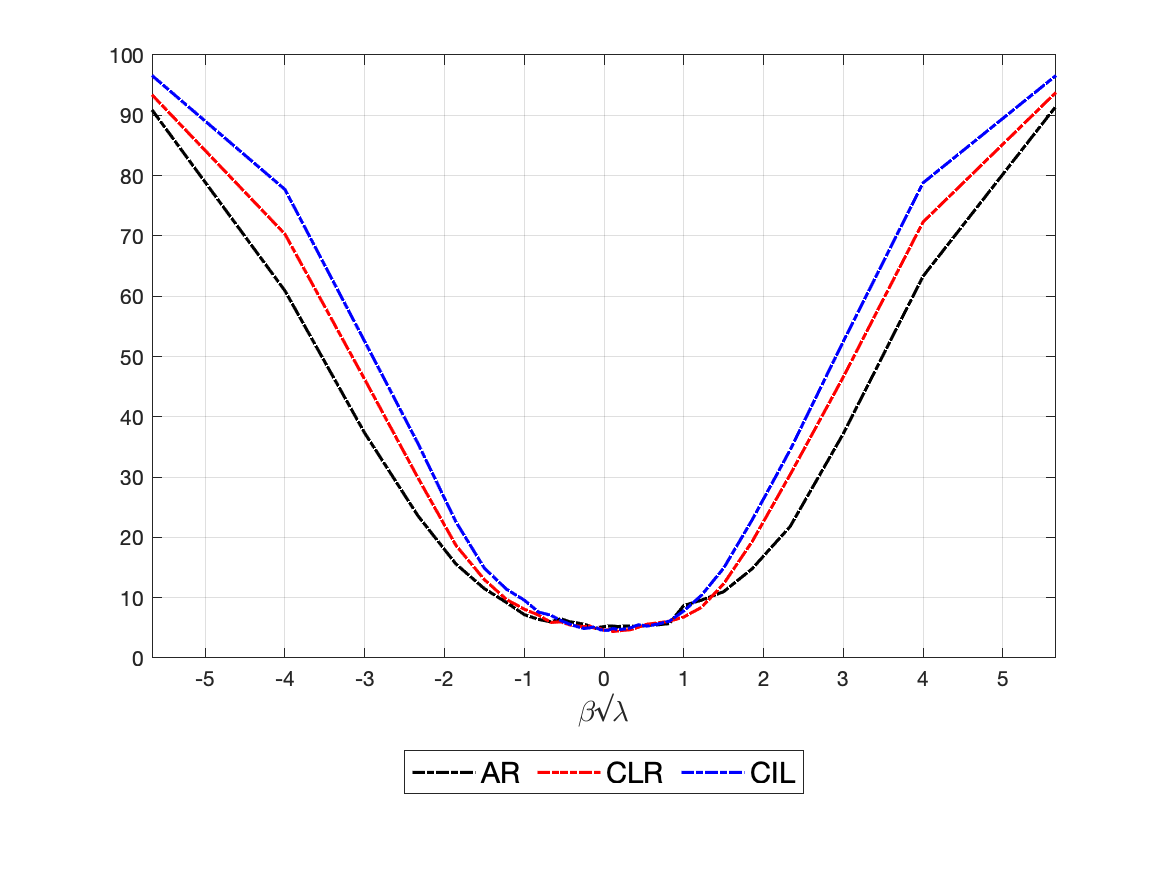}
\subcaption{NS with perturbation}
\end{subfigure}
\begin{subfigure}{.5\textwidth} \centering
\includegraphics[width=.9\linewidth]{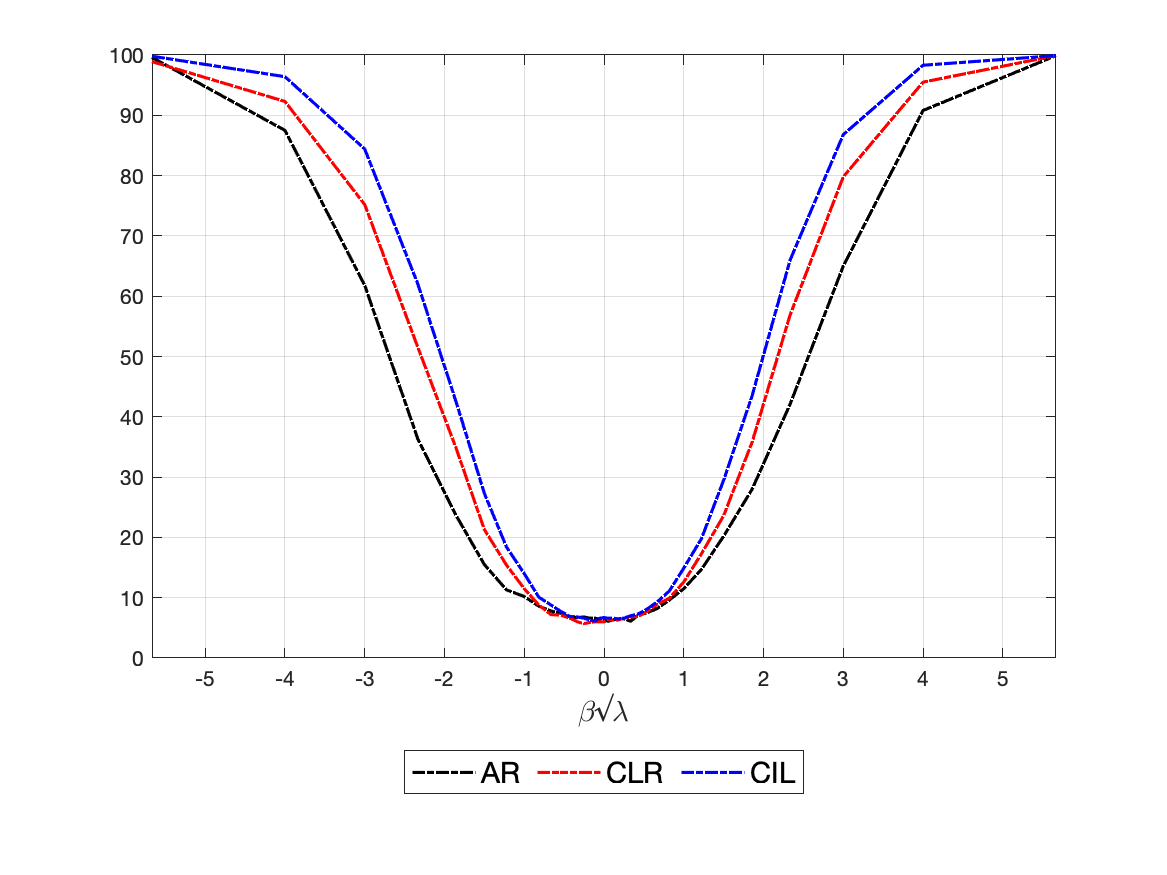}
\subcaption{GA Design}
\end{subfigure}
\begin{subfigure}{.5\textwidth}
\centering
\includegraphics[width=.9\linewidth]{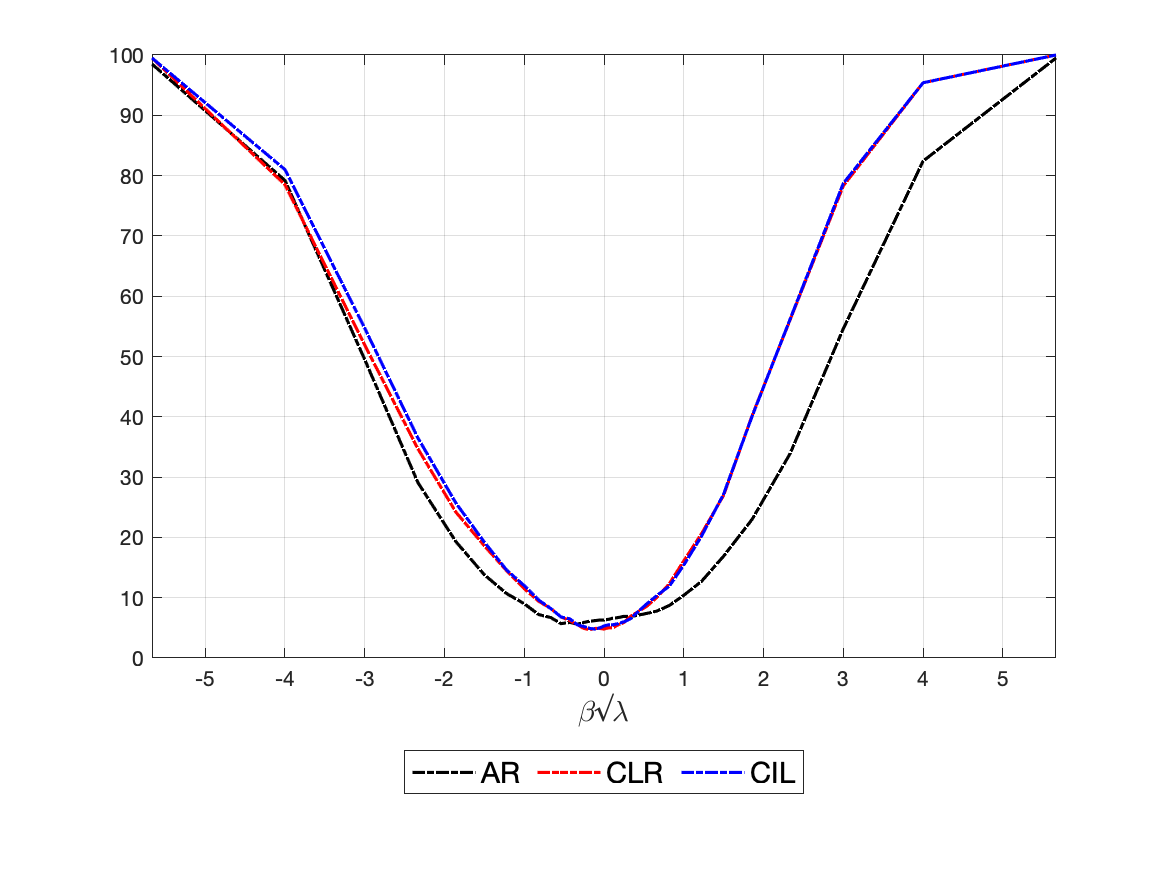}
\subcaption{NK Design}
\end{subfigure}
\end{figure}

We now briefly discuss power. More extensive power comparisons are reported
in the supplement (due to computational time for the CLR, these additional
comparisons use only 200 Monte Carlo replications for power and 200
simulations for conditional quantiles). Figure \ref{fig: HAC k=5} presents
power for the AR, CLR, CIL, and CIL$_{0}$ tests when $k=5$ and $\lambda /k=2$%
. We consider all four sets of simulations: NS design, NS design with a
perturbation, GA design, and NK design. As before, the CIL$_{0}$ test can be
biased, while the CLR and CIL\ tests dominate the AR test. In general, the
CIL\ test outperforms the CLR test. The power difference can be as large as
15\% for these specific designs. For example, the CIL\ test can have power
near 85\% when the CLR test rejects the null about 70\% of the time. This
difference happens even when we implement the infeasible version of the CLR
test which includes $\beta $ as one of the initial points.

The more technical sections of the paper are next. Section \ref{Sec
Kronecker} builds upon and connects with the work of AMS06. Section \ref{Sec
Invariant Tests} derives the CIL test.

\section{Kronecker Variance Matrix \label{Sec Kronecker}}

We first consider the special case where $\Sigma =\Omega \otimes \Phi $ with 
$\Omega $ a $2\times 2$ matrix and $\Phi $ a $k\times k$ matrix. The
Kronecker product framework is particularly interesting for two reasons.
First, we find the maximal invariant, taking into consideration a
transformation of $\Omega $ which is known but not fixed. This yields the
same data reduction from $S$ and $T$ as that obtained by AMS06 under the
assumption that $\Omega $ is known and fixed. This result is striking as 
\emph{the AMS06 approach does not hold} for general $\Sigma $, but \emph{%
ours does}. Second, AMS06 do not rule out the possibility that the test
depends on $\Omega $ beyond the statistics $S$ and $T$, because AMS06 treat $%
\Omega $ as being fixed. Our framework instead shows that invariant tests
should not depend on $\Omega $ \emph{at all}.

The $S$ and $T$ statistics in (\ref{(Defns of S and T)}) simplify to the
original statistics of \citet{Moreira02,Moreira09a} and AMS06 for the
homoskedastic model. When $\Sigma =\Omega \otimes \Phi $, the statistics $S$
and $T$ become%
\begin{eqnarray}
S &=&\Phi ^{-1/2}(Z^{\prime }Z)^{-1/2}Z^{\prime }Yb_{0}\cdot (b_{0}^{\prime
}\Omega b_{0})^{-1/2}\text{ and }  \label{(Defns of S and T kron)} \\
T &=&\Phi ^{-1/2}(Z^{\prime }Z)^{-1/2}Z^{\prime }Y\Omega ^{-1}a_{0}\cdot
(a_{0}^{\prime }\Omega ^{-1}a_{0})^{-1/2}.  \notag
\end{eqnarray}%
Their distribution is given by%
\begin{equation}
S\sim N\left( c_{\beta }\Phi ^{-1/2}\mu ,I_{k}\right) \text{ and }T\sim
N\left( d_{\beta }\Phi ^{-1/2}\mu ,I_{k}\right)  \label{(Dist S and T kron)}
\end{equation}%
with $c_{\beta }=(\beta -\beta _{0})\cdot (b_{0}^{\prime }\Omega
b_{0})^{-1/2}$ and $d_{\beta }=a^{\prime }\Omega ^{-1}a_{0}\cdot
(a_{0}^{\prime }\Omega ^{-1}a_{0})^{1/2}$. AMS06 develop the theory of
invariant tests by treating $\Omega $ as known and fixed. Even if $\Phi $ is
known, the parameter $\mu _{\Phi }=\Phi ^{-1/2}\mu $ is unknown, because $%
\mu $ is unknown. Hence, AMS06's invariance argument applies to the new
parameter $\mu _{\Phi }=\Phi ^{-1/2}\mu $. Specifically, let $h_{1}\in 
\mathcal{O}\left( k\right) $, the group of orthogonal matrices with matrix
multiplication as the group operator. The corresponding transformation in
the sample space is 
\begin{equation}
h_{1}\circ \left[ S:T\right] =h_{1}\cdot \left[ S:T\right] .
\end{equation}%
The associated transformation in the parameter space is%
\begin{equation}
h_{1}\circ \left( \beta ,\mu _{\Phi }\right) =\left( \beta ,h_{1}.\mu _{\Phi
}\right) .
\end{equation}%
The transformation does not change $\beta $, so our testing problem is
preserved. As argued before, this means that the test statistic should be an
invariant statistic (under the transformation $h_{1}$).

The maximal invariant statistic for the orthogonal transformation is%
\begin{equation}
Q=\left[ 
\begin{array}{cc}
Q_{S} & Q_{ST} \\ 
Q_{ST} & Q_{T}%
\end{array}%
\right] =\left[ 
\begin{array}{cc}
S^{\prime }S & S^{\prime }T \\ 
S^{\prime }T & T^{\prime }T%
\end{array}%
\right] .  \label{(Q def)}
\end{equation}%
That is, any invariant test depends on the data only through $Q$. The
density of $Q$ at $q$ for the parameters $\beta $ and $\lambda =\mu _{\Phi
}^{\prime }\mu _{\Phi }$ is given by 
\begin{eqnarray}
&&f_{\beta ,\lambda }(q_{S},q_{ST},q_{T})\overset{}{=}K_{0}\exp (-\lambda
(c_{\beta }^{2}+d_{\beta }^{2})/2)\left\vert q\right\vert ^{(k-3)/2} \\
&&\hspace{1.03in}\times \exp (-(q_{S}+q_{T})/2)(\lambda \xi _{\beta
}(q))^{-(k-2)/4}I_{(k-2)/2}(\sqrt{\lambda \xi _{\beta }(q)}),  \notag
\end{eqnarray}%
where $K_{0}^{-1}=2^{(k+2)/2}pi^{1/2}\Gamma _{(k-1)/2}$, $\Gamma _{(\cdot )}$
is the gamma function, $I_{(k-2)/2}(\cdot )$ denotes the modified Bessel
function of the first kind, and 
\begin{equation}
\xi _{\beta }(q)=c_{\beta }^{2}q_{S}+2c_{\beta }d_{\beta }q_{ST}+d_{\beta
}^{2}q_{T}.  \label{(Defn of Xi_Beta)}
\end{equation}

AMS06 further shows that another group, given by \emph{sign }%
transformations, preserves $H_{0}:\beta =\beta _{0}$ against $H_{0}:\beta
\neq \beta _{0}$. Consider the group $\mathcal{O}\left( 1\right) $, which
contains only two elements: $h_{2}\in $ $\left\{ -1,1\right\} $. For $%
h_{2}=-1$, the data transformation is given by%
\begin{equation}
h_{2}\circ \left[ S:T\right] =\left[ -S:T\right]
\end{equation}%
(by the definition of a group, the parameter remains unaltered at $h_{2}=1$%
). This yields a transformation in the maximal invariant space for $h_{2}$: 
\begin{equation}
h_{2}\circ \left( Q_{S},Q_{ST},Q_{T}\right) =\left(
Q_{S},h_{2}Q_{ST},Q_{T}\right) .
\end{equation}%
The maximal invariant for the joint transformation $h=\left(
h_{1},h_{2}\right) $ is the vector with components $Q_{S}$, $Q_{ST}^{2}$,
and $Q_{T}$. In principle, the tests can depend on $\Omega $ with
homoskedastic errors. As we will see in Theorem \ref{2sided Kronecker Thm},
we are able to eliminate the dependence on the variance as well, and show
the triad $Q_{S}$, $Q_{ST}^{2}$, and $Q_{T}$ is the maximal invariant for $%
g=\left( g_{1},g_{2}\right) $ in the case of known, but not fixed, variance.

\subsection{Instrument Transformation}

The orthogonal transformation argument of AMS06 is originally designed for
homoskedastic errors. For the general Kronecker case, both $\Phi ^{1/2}S$
and $\Phi ^{1/2}T$ (which are equivalent to the original statistics of
AMS06) have variance $\Phi $. Because their methodology assumes the variance
to be fixed, their orthogonal transformation would not work, in general,
because the variance would change. We could manually standardize their
statistics by $\Phi ^{-1/2}$ to obtain our statistics $S$ and $T$, and apply
the orthogonal group, as done earlier.\footnote{%
We could look instead at $g_{1}\in \mathcal{G}_{L}\left( k\right) $ such
that $g_{1}\Phi g_{1}^{\prime }=\Phi $. This yields $g_{1}=\Phi
^{1/2}h_{1}\Phi ^{-1/2}$. This is the same as transforming the data to $%
R_{\Phi }=\Phi ^{-1/2}R$, applying the orthogonal transformations, and
transforming the data back to $R$.} An alternative solution is to allow $%
\Phi $ to be known, but for it to change as we transform the data. For
example, take the special case in which $\Phi $ is a diagonal matrix. If we
were to permute the entries of $S$ and $T$ jointly, perhaps we should allow
the permutation of the diagonal entries of $\Phi $ as well. Formally, we
will take the variance $\Sigma =\Omega \otimes \Phi $ as part of both data
and parameter spaces.

For the special case in which $\Sigma =\Omega \otimes \Phi $, the
distribution of $R_{0}$ is given by 
\begin{equation}
R_{0}\sim N\left( \mu \left( \Delta ,1\right) ,\Omega _{0}\otimes \Phi
\right) ,
\end{equation}%
where $\Delta =\beta -\beta _{0}$, $\Sigma _{0}=\Omega _{0}\otimes \Phi $,
and $\Omega _{0}=B_{0}^{\prime }\Omega B_{0}$. The data are the realizations 
$(R_{0},\Omega _{0},\Phi )$ and the parameters are $(\Delta ,\mu ,\Omega
_{0},\Phi )$. The matrices $\Omega _{0},\Phi $ are assumed to be known, but
not fixed. Thus, $\Omega _{0},\Phi $ are both parameters and part of the
data, simultaneously.

We introduced the $g_{1}\in G_{L}\left( k\right) $ transformation in Section %
\ref{Sec Transformations}. Its action on the sample space is given by%
\begin{equation}
g_{1}\circ \left( R_{0},\Omega _{0},\Phi \right) =\left( g_{1}R_{0},\Omega
_{0},g_{1}\Phi g_{1}^{\prime }\right) .
\label{(left action sample-space Kronecker)}
\end{equation}%
We note that%
\begin{equation}
g_{1}R_{0}\sim N\left( g_{1}\mu \left( \Delta ,1\right) ,\Omega _{0}\otimes
g_{1}\Phi g_{1}^{\prime }\right) ,
\end{equation}%
so the corresponding action on the parameter space is%
\begin{equation}
g_{1}\circ \left( \Delta ,\mu ,\Omega _{0},\Phi \right) =\left( \Delta
,g_{1}\mu ,\Omega _{0},g_{1}\Phi g_{1}^{\prime }\right) .
\label{(left action parameter-space Kronecker)}
\end{equation}%
We now show that the matrix

\begin{equation}
Q=\left[ S:T\right] ^{\prime }\left[ S:T\right] =\left[ 
\begin{array}{cc}
S^{\prime }S & S^{\prime }T \\ 
S^{\prime }T & T^{\prime }T%
\end{array}%
\right] ,
\end{equation}%
together with $\Omega _{0}$ itself, is the maximal invariant statistic. That
is, any other invariant statistic can be written as a function of $\left(
Q,\Omega _{0}\right) $. The distribution of the maximal invariant depends
only on the concentration parameter $\lambda $, the parameter of interest $%
\beta $, and $\Omega _{0}$ itself.

\bigskip

\begin{theorem}
\label{Kronecker Thm} For the group actions in (\ref{(left action
sample-space Kronecker)}) and (\ref{(left action parameter-space Kronecker)}%
):\newline
\emph{(i)} The maximal invariant in the sample space is given by $\left(
Q,\Omega _{0}\right) $; and\newline
\emph{(ii) }The maximal invariant in the parameter space is given by $\left(
c_{\beta }^{2}\lambda ,c_{\beta }d_{\beta }\lambda ,d_{\beta }^{2}\lambda
,\Omega _{0}\right) $.
\end{theorem}

\textbf{Comments: 1. }The data $\left( \left[ S:T\right] ,\Omega _{0},\Phi
\right) $ is a one-to-one transformation from the primitive data $\left(
R,\Omega ,\Phi \right) $. Hence, there is no loss of generality in using the 
\emph{pivotal} statistic $S$ and the \emph{complete} statistic $T$ instead
of using $R$ (or $R_{0}$).

\textbf{2.} There is a one-to-one mapping between $\Omega _{0}$ and $\Omega $%
. Hence, $\left( Q,\Omega \right) $ is a maximal invariant as well. We
continue to use $\Omega _{0}$ because it is useful to find a maximal
invariant for the two-sided transformations to be considered next.

\textbf{3.} The statistic $Q$ is the maximal invariant based on the compact
orthogonal group on $\left[ S:T\right] $, which is a straightforward
application of AMS06. We instead allow the much larger, noncompact group of
nonsingular matrices with unitary determinant. The data also contain the
variance components given by $\Omega _{0}$ and $\Phi $. Because the group $%
\mathcal{G}_{L}\left( k\right) $ is not \emph{amenable}, the Hunt-Stein
theorem is not applicable, and we do not necessarily obtain a minimax
result. This is in contrast to \citet{Chamberlain07}, who builds on the fact
that the orthogonal group is compact.

\textbf{4.} The component $\Phi $ completely vanishes as the noncompact
group $\mathcal{G}_{L}\left( k\right) $ acts \emph{transitively} on $\Phi $.
Hence, the matrix $\Phi $ is not part of the maximal invariant.

\subsection{Two-Sided Transformation \label{Sec two-sided}}

We now apply the $g_{2}\in \mathcal{G}_{T}\left( 2\right) $ transformation
introduced in Section \ref{Sec Transformations}. The two-sided
transformation in the Kronecker model is given by 
\begin{equation}
g_{2}\circ \left( R_{0},\Omega _{0},\Phi \right) =\left( R_{0}g_{2}^{\prime
},g_{2}\Omega _{0}g_{2}^{\prime },\Phi \right) ,
\label{(right action sample-space Kronecker)}
\end{equation}%
where $g_{2}\in \mathcal{G}_{T}\left( 2\right) $, the group of nonsingular
lower triangular $2\times 2$ matrices. The transformation in the parameter
space is 
\begin{equation}
g\circ \left( \Delta ,\mu ,\Sigma _{0}\right) =\left( \frac{\Delta g_{11}}{%
\Delta g_{21}+g_{22}},g_{1}\mu \left( \Delta g_{21}+g_{22}\right) ,\left(
g_{2}\otimes g_{1}\right) \Sigma _{0}\left( g_{2}^{\prime }\otimes
g_{1}^{\prime }\right) \right) .
\label{(right action
parameter-space Kronecker)}
\end{equation}%
Theorem \ref{2sided Kronecker Thm} finds the maximal invariant based on $%
g_{1}\in \mathcal{G}_{L}\left( k\right) $ and $g_{2}\in \mathcal{G}%
_{T}\left( 2\right) $.

\bigskip

\begin{theorem}
\label{2sided Kronecker Thm} For the data group actions defined in (\ref%
{(left action sample-space Kronecker)}) and (\ref{(right action sample-space
Kronecker)}), and the parameter actions in (\ref{(left action
parameter-space Kronecker)}) and (\ref{(right action parameter-space
Kronecker)}), we find\newline
\emph{(i)} The induced group action of $g_{2}$ on the space $\left( \left[
S:T\right] ,\Omega _{0},\Phi \right) $ is%
\begin{equation*}
g_{2}\circ \left( \left[ S:T\right] ,\Omega _{0},\Phi \right) =\left( \left[
sgn\left( g_{11}\right) S:sgn\left( g_{22}\right) T\right] ,g_{2}\Omega
_{0}g_{2}^{\prime },\Phi \right) ;
\end{equation*}%
\emph{(ii)} The data maximal invariant to $g=\left( g_{1},g_{2}\right) $ is 
\begin{equation*}
\left( Q_{S},Q_{T},Q_{ST}^{2}\right) ;
\end{equation*}%
\emph{(iii)} The induced group action by $g_{2}$ on the parameter functions $%
\left( c_{\beta },d_{\beta },\mu ,\Omega _{0},\Phi \right) $ is given by%
\begin{equation*}
g_{2}\circ \left( c_{\beta }\mu ,d_{\beta }\mu ,\Omega _{0},\Phi \right)
=\left( sgn\left( g_{11}\right) c_{\beta }\mu ,sgn\left( g_{22}\right)
d_{\beta }\mu ,g_{2}\Omega _{0}g_{2}^{\prime },\Phi \right) ;\text{ and}
\end{equation*}%
\emph{(iv)} The parameter maximal invariant to $g=\left( g_{1},g_{2}\right) $
is%
\begin{equation*}
\left( c_{\beta }^{2}\lambda ,d_{\beta }^{2}\lambda ,\left\vert c_{\beta
}d_{\beta }\right\vert \lambda \right) .
\end{equation*}
\end{theorem}

\textbf{Comments: 1. }The parameters $\beta $ and $\Omega $ remain unchanged
by the action (\ref{(left action parameter-space Kronecker)}). Because the
parameters $c_{\beta }$ and $d_{\beta }$ depend only on $\beta $ and $\Omega 
$, they are preserved as well. The result now follows trivially because $%
g_{1}\circ \left( \mu ,\Omega ,\Phi \right) =\left( g_{1}\mu ,\Omega
,g_{1}\Phi g_{1}^{\prime }\right) $.

\textbf{2.} We note that $g_{21}$ may be different from zero. Hence, the
group of transformations is larger than scale multiplication to each entry
of the vector $\left( \Delta ,1\right) $. A naive generalization for the
sign group of transformations by AMS06 to our setup is a diagonal matrix $%
g_{2}$. In the online appendix, we show that some invariant tests based on
the associated maximal invariant can behave as one-sided tests. Hence, we
illustrate the importance of finding the largest group of transformations
before deriving invariant tests.

\bigskip

These actions are defined using the reduced-form matrix $\Omega $. For the
homoskedastic model, we could analyze the transformations in the
structural-form matrix%
\begin{equation}
\Psi =\left[ 
\begin{array}{cc}
\sigma _{uu} & \sigma _{u2} \\ 
\sigma _{u2} & \sigma _{22}%
\end{array}%
\right] .
\end{equation}%
One may wonder if there are actually symmetries in the original model. This
turns out to be true. In fact, the action in the structural-form variance
matrix has a very simple structure.

\bigskip

\begin{proposition}
\label{Structural Prop} The group action on the reduced-form matrix $\Omega $
induces an action on the structural-form matrix $\Psi $:%
\begin{eqnarray*}
g_{2}\circ \left( \Delta ,\lambda ,\Psi \right) &=&\left( \frac{\Delta g_{11}%
}{\Delta g_{21}+g_{22}},\left( \Delta g_{21}+g_{22}\right) ^{2}\lambda
,\Gamma \Psi \Gamma ^{\prime }\right) \text{, where} \\
\Gamma &=&\left[ 
\begin{array}{cc}
\left( \Delta g_{21}+g_{22}\right) ^{-1}g_{11}g_{22} & 0 \\ 
g_{21} & \Delta g_{21}+g_{22}%
\end{array}%
\right] .
\end{eqnarray*}
\end{proposition}

\textbf{Comment: }Take $\beta _{0}=0$. When $g_{11}=-1,$ $g_{21}=0,$ and $%
g_{22}=1$, we have $g_{2}\circ \left( v_{1},v_{2}\right) =\left(
-v_{1},v_{2}\right) $. Therefore, $\sigma _{11}$ and $\sigma _{22}$ are
preserved while $\sigma _{12}$ changes sign. Since $\sigma _{12}=\sigma
_{u2}+\sigma _{22}\beta $, the new value for the structural-form covariance
scalar, $-\sigma _{u2}$, and the new value of the parameter, $-\beta $,
comprise the only transformation that works for any value of $\sigma _{22}$.

\bigskip

A corollary of our theory is that the AR test is UMPI when structural-form
variance is fixed. This optimality result is novel and important. All
optimality theorems for the AR test, so far, assume the reduced-form
variance to be fixed (\citet{Moreira02,Moreira09a}, AMS06, and %
\citet{MoreiraMoreira19}).

\section{Invariant Tests \label{Sec Invariant Tests}}

We now use the group of transformations by $g=\left( g_{1},g_{2}\right) $ to
develop the CIL\ test. Recall that the data consist of $R_{0}$ and $\Sigma
_{0}$, where $R_{0}$ has a normal distribution and the distribution of $%
\Sigma _{0}$ is degenerate. So, the density of the data is the product of
two parts. The first part is the normal distribution of $R_{0}$, which is
absolutely continuous with respect to the Lebesgue measure. The second part
is the degenerate distribution of $\Sigma _{0}$, which is absolutely
continuous with respect to the counting measure. Understanding how the
density changes with the data transformation is important for the
development of the CIL\ test.

The density of $R_{0}=\left[ R_{1}:R_{2}\right] $ evaluated at $r_{0}=\left[
r_{1}:r_{2}\right] $ is given by%
\begin{equation}
f_{R}\left( r_{0};\Delta ,\mu ,\Sigma _{0}\right) =\left( 2pi\right)
^{-k}\left\vert \Sigma _{0}\right\vert ^{-1/2}\exp \left\{ -\frac{1}{2}\left[
\begin{array}{c}
r_{1}-\mu \Delta \\ 
r_{2}-\mu%
\end{array}%
\right] ^{\prime }\Sigma _{0}^{-1}\left[ 
\begin{array}{c}
r_{1}-\mu \Delta \\ 
r_{2}-\mu%
\end{array}%
\right] \right\} .
\end{equation}

As in Theorem \ref{2sided Kronecker Thm}, we consider the groups of
instrument transformations $g_{1}$ and two-sided transformations $g_{2}$
together, so that we have the joint transformation $g=\left(
g_{1},g_{2}\right) $ defined in Section \ref{Sec Transformations}, where $%
g_{1}\in \mathcal{G}_{L}\left( k\right) $ and $g_{2}\in \mathcal{G}%
_{T}\left( 2\right) $, and the associated transformation $g\circ \left(
\Delta ,\mu ,\Sigma _{0}\right) $ in the parameter space.

Basic algebraic manipulations show that%
\begin{equation}
f_{R}\left( g\circ r_{0};g\circ \left( \Delta ,\mu ,\Sigma _{0}\right)
\right) =f_{R}\left( \left[ r_{1}:r_{2}\right] ;\Delta ,\mu ,\Sigma
_{0}\right) \left\vert g_{2}\right\vert ^{-k}\left\vert g_{1}\right\vert
^{-2}
\end{equation}%
because 
\begin{equation}
\left\vert \left( g_{2}\otimes g_{1}\right) \Sigma _{0}\left( g_{2}^{\prime
}\otimes g_{1}^{\prime }\right) \right\vert =|g_{2}|^{2k}|g_{1}|^{4}|\Sigma
_{0}|.
\end{equation}%
Therefore,%
\begin{equation}
f_{R}\left( r_{0};\Delta ,\mu ,\Sigma _{0}\right) =f_{R}\left( g\circ
r_{0};g\circ \left( \Delta ,\mu ,\Sigma _{0}\right) \right) \chi _{0}\left(
g\right) ,
\end{equation}%
where $\chi _{0}\left( g\right) =\chi _{1}\left( g_{1}\right) \chi
_{2}\left( g_{2}\right) $ for the sub-group multipliers $\chi _{1}\left(
g_{1}\right) =\left\vert g_{1}\right\vert ^{2}$ and $\chi _{2}\left(
g_{2}\right) =\left\vert g_{2}\right\vert ^{k}$. Hence, the density of $%
R_{0} $ is relatively invariant with multiplier $\chi _{0}(g)$.

Of course, the action $g\in \mathcal{G}_{L}\left( k\right) \times \mathcal{G}%
_{T}\left( 2\right) $ is not \emph{proper}.\footnote{%
See Definition 5.1 of \citet{Eaton89} for a formal statement on a group
acting properly on the sample space. In our case, it is trivial that the
action by $g$ is not proper, since we can multiply $g_{1}$ and divide $g_{2}$
by the same constant.} We can impose $\left\vert g_{1}\right\vert =1$ so
that $\chi _{1}\left( g_{1}\right) =1$. In this case, $g_{1}\in \mathcal{S}%
_{L}\left( k\right) $, the group of invertible matrices with determinant
equal to one. Alternatively, we can use another standardization such as $%
g_{22}=1$. To develop the integrated likelihood invariant test, we use Haar
measures to obtain invariant tests. It is harder to work with the Haar
measure for $\mathcal{S}_{L}\left( k\right) $ than for $\mathcal{G}%
_{L}\left( k\right) $; see \citet{Dedic90}. On the other hand, it is
relatively simple to derive the Haar measure for $2\times 2$ lower
triangular matrices with $g_{22}=1$. For this reason, we prefer to impose a
restriction on $\mathcal{G}_{T}\left( 2\right) $.

For the second part, the data $\Sigma _{0}$ have a distribution that assigns
probability one to the value $\Sigma _{0}$ itself. Therefore, the density at
some arbitrary matrix value $\sigma _{0}$ is

\begin{equation}
f_{\Sigma }\left( \sigma _{0};\Sigma _{0}\right) =P_{\Sigma _{0}}\left(
\sigma _{0}=\Sigma _{0}\right) =I\left( \sigma _{0}=\Sigma _{0}\right) .
\label{(degenerate density eq)}
\end{equation}%
Using (\ref{(degenerate density eq)}), we have%
\begin{equation}
f_{\Sigma }(g\circ \sigma _{0};g\circ \Sigma _{0})=f_{\Sigma }\left( \left(
g_{2}\otimes g_{1}\right) \sigma _{0}\left( g_{2}^{\prime }\otimes
g_{1}^{\prime }\right) ;\left( g_{2}\otimes g_{1}\right) \Sigma _{0}\left(
g_{2}^{\prime }\otimes g_{1}^{\prime }\right) \right) =f_{\Sigma }\left(
\sigma _{0};\Sigma _{0}\right)
\end{equation}%
so that this density is invariant with multiplier 1.

The joint likelihood is then given by%
\begin{equation}
f\left( r_{0},\sigma _{0};\Delta ,\mu ,\Sigma _{0}\right) =f_{R}\left(
r_{0};\Delta ,\mu ,\Sigma _{0}\right) \cdot f_{\Sigma }\left( \sigma
_{0};\Sigma _{0}\right) ,
\end{equation}%
so that%
\begin{equation}
f\left( r_{0},\sigma _{0};\Delta ,\mu ,\Sigma _{0}\right) =f\left( g\circ
\left( r_{0},\sigma _{0}\right) ;g\circ \left( \Delta ,\mu ,\Sigma
_{0}\right) \right) \cdot \chi _{0}\left( g\right) ,
\label{invariance-model}
\end{equation}%
i.e. the likelihood is relatively invariant with multiplier $\chi _{0}(g)$.
Because the Lebesgue measure is relatively left invariant for the group $g$
with multiplier $\chi _{0}(g)$, the (relative) invariance of the likelihood
follows directly.

We use the invariance of the likelihood to propose a conditional weighted
likelihood ratio test. We also show that the AR, LM, CQLR, and CLR tests are
also invariant.

\subsection{Optimal Tests \label{Subsec Optimal Tests}}

Our goal in this section is to find optimal tests. Specifically, a test is
defined to be a measurable function $\phi \left( r_{0},\sigma _{0}\right) $
that is bounded by $0$ and $1$. For a given outcome, the test rejects the
null with probability $\phi \left( r_{0},\sigma _{0}\right) $ and accepts
the null with probability $1-\phi \left( r_{0},\sigma _{0}\right) $, e.g.,
the Anderson-Rubin test is simply $I\left( AR>c\left( k\right) \right) $
where $I\left( \cdot \right) $ is the indicator function. The test is said
to be nonrandomized if $\phi $ takes only values $0$ and $1$; otherwise, it
is called a randomized test. The rejection probability is given by%
\begin{equation}
E_{\Delta ,\mu ,\Sigma _{0}}\phi \left( R_{0},\Sigma _{0}\right) \equiv \int
\phi \left( r_{0},\sigma _{0}\right) f\left( r_{0},\sigma _{0};\Delta ,\mu
,\Sigma _{0}\right) \text{ }dr_{0}\text{ }\eta \left( d\sigma _{0}\right) ,
\label{(power fn eq)}
\end{equation}%
where $\eta $ is the counting measure. The rejection probability (\ref%
{(power fn eq)}) simplifies to%
\begin{eqnarray}
E_{\Delta ,\mu ,\Sigma _{0}}\phi \left( R_{0},\Sigma _{0}\right) &=&\int
\phi \left( r_{0},\sigma _{0}\right) f_{R}\left( r_{0};\Delta ,\mu ,\Sigma
_{0}\right) f_{\Sigma }\left( \sigma _{0};\Sigma _{0}\right) \text{ }dr_{0}%
\text{ }\eta \left( d\sigma _{0}\right)  \notag \\
&=&\int \phi \left( r_{0},\Sigma _{0}\right) f_{R}\left( r_{0};\Delta ,\mu
,\Sigma _{0}\right) \text{ }dr_{0}.
\end{eqnarray}

The rejection probability $E_{\Delta ,\mu ,\Sigma _{0}}\phi \left(
R_{0},\Sigma _{0}\right) $ taken as a function of $\Delta $, $\mu $, and $%
\Sigma _{0}$ gives the power curve for the test $\phi $. In particular, $%
E_{0,\mu ,\Sigma _{0}}\phi \left( R_{0},\Sigma _{0}\right) $ gives the null
rejection probability.

Let the parameter space for $\Delta ,\mu ,\sigma _{0}$ be denoted by $\Theta 
$, with $\sigma $-field the intersection of $\Theta $ and sets in $\mathcal{B%
}^{k+1}\times \{\Sigma _{0}\}$. Let $w$ be a measure on that $\sigma $%
-field. We average the power curve over the parameter space to obtain the
weighted average power with weights that are given by the measure $w$. By
Tonelli's theorem, the weighted average power is 
\begin{equation}
E_{w}\phi \left( R_{0},\Sigma _{0}\right) =\int E_{\Delta ,\mu ,\Sigma
_{0}}\phi \left( R_{0},\Sigma _{0}\right) dw\left( \Delta ,\mu ,\Sigma
_{0}\right) .
\end{equation}%
If the weights are such that for $B\times \{\Sigma _{0}\}$, 
\begin{equation}
w\left( B\times \{\Sigma _{0}\}\right) =w_{R}\left( B\right) \cdot w_{\Sigma
}\left( \{\sigma _{0}\}\right) ,
\end{equation}%
where $B\in \mathcal{B}^{k+1}$ and $w_{\Sigma }\left( \{\sigma _{0}\}\right) 
$ has unitary mass on $\{\Sigma _{0}\}$, then%
\begin{equation}
E_{w}\phi \left( R_{0},\Sigma _{0}\right) =\int \phi \left( r_{0},\Sigma
_{0}\right) f_{w_{R}}\left( r_{0},\Sigma _{0}\right) \text{ }dr_{0},
\end{equation}%
where $f_{w_{R}}\left( r_{0},\Sigma _{0}\right) $ is defined as 
\begin{equation}
f_{w_{R}}\left( r_{0},\Sigma _{0}\right) =\int f_{R}\left( r_{0};\Delta ,\mu
,\Sigma _{0}\right) \text{ }dw_{R}\left( \Delta ,\mu \right) .
\end{equation}

For a given weight $w$, we seek optimal similar tests 
\begin{equation}
\max_{0\leq \phi \leq 1}E_{w}\phi \left( R_{0},\Sigma _{0}\right) \text{,
where }E_{0,\mu ,\Sigma _{0}}\phi \left( R_{0},\Sigma _{0}\right) =\alpha
,\forall \mu .  \label{(WAP similar)}
\end{equation}%
The next proposition finds the WAP test.

\bigskip

\begin{proposition}
\label{WAP similar Prop} The optimal test in (\ref{(WAP similar)}) rejects
the null when%
\begin{equation}
\frac{f_{w_{R}}\left( r_{0},\Sigma _{0}\right) }{f_{S}\left( s\right) }%
>\kappa \left( t,\Sigma _{0}\right) ,  \label{(WAP similar hT)}
\end{equation}%
where $f_{S}\left( s\right) =\left( 2pi\right) ^{-k/2}e^{-s^{\prime }s/2}$
is the density of the statistic $S$ under the null.
\end{proposition}

\textbf{Comment: }Because $T$ is sufficient for $\mu $ under the null, we
condition on $T=t$. The dependence of the test statistic on $t$ is absorbed
in the critical value of the test.

\bigskip

For arbitrary weights, the WAP similar test is not guaranteed to have
overall good power in finite samples. In particular, the power can be near
zero for parts of the parameter space (as happens with the CIL$_{0}$ test
for $k>2$). We circumvent this problem by carefully choosing weights $w$ so
that the test given by (\ref{(WAP similar hT)}) is invariant. The CIL\ test
behaves as a two-sided test, and so, it does not suffer the criticism by %
\citet{MoreiraMoreira19}.

\subsection{Similar Invariant Tests \label{Subsec The Statistic}}

Invariance of conditional tests follows from the relative invariance of test
statistics.

\bigskip

\begin{definition}
\label{Left Invariant Stat Def} A statistic $\psi $ is relatively (left)
invariant to $g$ with multiplier $\chi $ if%
\begin{equation*}
\psi \left( g\circ \left( s,t,\sigma _{0}\right) \right) =\chi \left(
g\right) \cdot \psi \left( s,t,\sigma _{0}\right) ,
\end{equation*}%
for any $\left( s,t,\sigma _{0}\right) $.
\end{definition}

\bigskip

Proposition \ref{Invariance of Conditional Tests Prop} establishes the
invariance of the conditional test if the test statistic is relatively
invariant.

\bigskip

\begin{proposition}
\label{Invariance of Conditional Tests Prop} Suppose that $\psi \left(
S,t,\Sigma _{0}\right) $ is a continuous random variable under $H_{0}:\Delta
=0$ for every $t$. Define $\kappa _{\psi }\left( t,\Sigma _{0}\right) $ to
be the $1-\alpha $ quantile of the null distribution of $\psi \left(
S,t,\Sigma _{0}\right) $. Then the following hold:\newline
\emph{(i)} The conditional test $\phi \left( s,t,\Sigma _{0}\right) $ that
rejects the null when%
\begin{equation*}
\psi \left( s,t,\Sigma _{0}\right) >\kappa _{\psi }\left( t,\Sigma
_{0}\right)
\end{equation*}%
is similar at level $\alpha $;\newline
\emph{(ii)} If $\psi \left( g\circ \left( s,t,\Sigma _{0}\right) \right) $
is relatively invariant under $g\in \mathcal{G}_{L}\left( k\right) \times
G_{T}\left( 2\right) $ with multiplier $\chi $, then $\kappa _{\psi }\left(
t,\Sigma _{0}\right) $ is itself relatively invariant with multiplier $\chi $%
; and\newline
(iii) The conditional test $\phi \left( s,t,\Sigma _{0}\right) $ is
invariant.
\end{proposition}

\textbf{Comments: 1.} Careful examination of the proof shows that invariance
of the conditional quantile does not depend on the group transformation
used. It is also applicable to other models as long as there is a sufficient
statistic, e.g. here under the null, that is boundedly complete.

\textbf{2.} The comment above explains why the conditional quantile of the $%
LR$ statistic depends only on $T^{\prime }T$ in the homoskedastic case. The
LR statistic does not depend on $\Omega _{0}$ at all, and $T^{\prime }T$ is
the maximal invariant to orthogonal transformations $h_{1}\circ T=h_{1}\cdot
T$. This is consistent with the results of \citet{Moreira03} and AMS06, but
with no need to use pivotal statistics and independence.

\bigskip

Before showing that the CIL\ test is invariant and is the limit of
conditional WAP tests, as given by (\ref{(WAP similar hT)}), we establish
that the $AR$, $LM$, $LR$, and $QLR$ statistics are invariant.

\begin{proposition}
\label{Invariant Test Stats Prop} The $AR$, $LM$, $LR$, and $QLR$ statistics
are invariant to $g=\left( g_{1},g_{2}\right) \in \mathcal{G}_{L}\left(
k\right) \times \mathcal{G}_{T}\left( 2\right) $.
\end{proposition}

\textbf{Comment:} Close inspection shows the proof of invariance of the $LR$
statistic is very general. It works for any model in the presence of
symmetries which preserve the testing problem.

\subsection{An Invariant WAP Similar Test \label{Subsec Invariant WAP}}

The goal is to obtain a WAP invariant similar test in the over-identified
model ($k>1$). This entails finding weights so that the final test is
relatively invariant.

\bigskip

\begin{definition}
\label{Relatively Invariant Def} A measure $m$ is relatively (left)
invariant with multiplier $\chi $ if%
\begin{equation*}
\int F\left( g^{-1}\circ \theta \right) m\left( d\theta \right) =\chi \left(
g\right) \int F\left( \theta \right) m\left( d\theta \right)
\end{equation*}%
for any real-valued continuous function $F$ with bounded support.
\end{definition}

\bigskip

We could apply this result for $\theta =\left( \Delta ,\mu ,\Sigma
_{0}\right) $. However, the parameter $\Sigma _{0}$ is known, but changes
according to the data transformation. Therefore, it is enough to allow $%
\theta $ to be the parameters $\left( \Delta ,\mu \right) $ only.

\bigskip

\begin{lemma}
\label{Relatively Invariant Lemma} The product measure $\left\vert \Delta
\right\vert ^{k-2}d\Delta \times $ $d\mu $ is relatively (left) invariant to 
$g=\left( g_{1},g_{2}\right) \in \mathcal{G}_{L}\left( k\right) \times
G_{T}\left( 2\right) $ with multiplier $\left\vert g_{1}\right\vert \cdot
\left\vert g_{11}\right\vert ^{k-1}$.
\end{lemma}

\bigskip

The next proposition shows that the conditional test is invariant and can be
evaluated with a single (and not multiple) integral.

\bigskip

\begin{theorem}
\label{WAP Invariant Cond Test Thm} The conditional test based on the test
statistic%
\begin{eqnarray}
IL &=&\int e^{-\frac{1}{2}\left[ vec\left( R_{0}\right) ^{\prime }\Sigma
_{0}^{-1/2}N_{\Sigma _{0}^{-1/2}(a_{\Delta }\otimes I_{k})}\Sigma
_{0}^{-1/2}vec\left( R_{0}\right) -T^{\prime }T\right] }  \label{WAP} \\
&&\times \left\vert \left( a_{\Delta }^{\prime }\otimes I_{k}\right) \Sigma
_{0}^{-1}\left( a_{\Delta }\otimes I_{k}\right) \right\vert
^{-1/2}\left\vert \Delta \right\vert ^{k-2}d\Delta  \notag
\end{eqnarray}%
is invariant and is the limit of a sequence of WAP tests defined in (\ref%
{(WAP similar hT)}).
\end{theorem}

\bigskip

In separate work, we address admissibility of the CIL test. Showing
admissibility based on invariant weights for non-amenable groups is done on
a case-by-case basis. This issue is analogous to that encountered for the
commonly-accepted and widely-used Hotelling $T^{2}$ statistic for testing
means of different populations. \citet{Stein55} addresses the admissibility
of the Hotelling $T^{2}$ statistic. This test relies on the same
non-amenable $\mathcal{G}_{L}\left( k\right) $ group considered here for the
HAC-IV model.

For the construction of the $IL$ statistic, both priors for $\theta $ and $%
\mu $ are improper. In the spirit of Theorem \ref{WAP Invariant Cond Test
Thm}, we need to consider sequences of weights for the alternative
hypothesis instead. For the nuisance parameter $\mu $, %
\citet{MoreiraMoreira19} and \citet{AndrewsMikusheva20} allow an
\textquotedblleft identification\textquotedblright\ parameter go to
infinity, so that the prior converges weakly to the Lebesgue measure. The
completeness theorem shows that any admissible test is the limit of Bayes
tests (sufficiency). However, is the limit of any sequence of Bayes tests
admissible (necessity)? \citet{Farrell68b,Farrell68a} considers a more
concrete version of Stein's proof of admissibility. For example, %
\citet{MoreiraMoreira13,MoreiraMoreira19} rely on Farrell's approach by
using subsequence arguments for admissibility of WAP similar tests.

\section{Conclusion and Extensions \label{Sec Conclusion}}

This paper shows the importance of distinguishing between parameters being
known or being fixed when showing the presence of symmetries in the HAC-IV
model. However, this distinction is applicable to many other models. The
existence of symmetries could be useful, as they could simplify inference
(e.g., data reduction by invariance) or enable us to find better estimators
and tests.

Econometricians are often interested in some parameters in the presence of
others. It is well-understood that knowing the value of a nuisance parameter
typically yields more efficient estimators and tests than estimating it. In
some cases, however, knowing or estimating the nuisance parameter yields the
same asymptotic efficiency. This feature can happen in parametric models in
cross-section or time-series data, as well as in semi-parametric models,
among others.

\begin{itemize}
\item Consider a linear regression when the error variance is unknown (up to
a parameter of fixed dimension). The generalized least-squares (GLS)
estimator is the optimal linear unbiased estimator. It enjoys asymptotic
efficiency among regular estimators. The feasible generalized least-squares
(FGLS) estimator is asymptotically efficient when we consistently estimate
the parametric error variance.

\item Take the predictive regression model where the explanatory variable
can be nearly integrated of order one. The asymptotic behavior of several
tests is the same whether the long-run variance matrix of the errors is
known or consistently estimated.

\item In the GMM model, we can consistently estimate the optimal weighting
matrix using a HAC estimator. Assuming the variance is known or estimated,
the GMM\ estimators are asymptotically equivalent and efficient.
\end{itemize}

These examples illustrate the caveats of estimating or testing by assuming
some parameters are \emph{known}. This natural simplification ironically
leads to complications when parameters are assumed to be \emph{fixed}. In
particular, it leads to the incorrect folk theorem that many models do not
present natural symmetries. Once we distinguish between the assumptions of
known versus fixed parameters, model symmetries can exist, contrary to
popular belief. We hope this new methodology will lead to new inferential
methods to apply to important econometric models.

\bibliographystyle{econometrica}
\bibliography{References}

\end{document}